\documentclass[9pt]{style1}

\newtheorem{theorem}{Theorem}[section]
\newtheorem{lemma}[theorem]{Lemma}
\newtheorem{corollary}[theorem]{Corollary}
\newtheorem{proposition}[theorem]{Proposition}

\theoremstyle{definition}
\newtheorem{definition}[theorem]{Definition}

\theoremstyle{remark}
\newtheorem{remark}[theorem]{Remark}

\def\tM{\widetilde M}

\def\sigmab{\bar{\sigma}}
\def\texp{\widetilde \exp}
\def\heat{\lf(\frac{\p}{\p t}-\Delta\ri)}
\def \b {\beta}

\def\lf{\left}
\def\ri{\right}

\def\a{\alpha}

\def\e{\epsilon}
\def\p{\partial}

\def\R{\Bbb R}

\def \D {\Delta}

\numberwithin{equation}{section}

\begin{document}

\title{Mean value theorems on manifolds}

\author{Lei Ni}
\address{Department of Mathematics, University of California at San Diego, La Jolla, CA 92093}

\email{lni@math.ucsd.edu}
\thanks{The author was supported in part by NSF Grants  and
an Alfred P. Sloan Fellowship, USA}


\date{December 2005}


\begin{abstract} We derive  several mean value formulae on
manifolds, generalizing the classical one for harmonic functions on
Euclidean spaces as well as later results of Schoen-Yau,
Michael-Simon, etc, on curved Riemannian manifolds. For the heat
equation a mean value theorem with respect to `heat spheres' is
proved for heat equation with respect to evolving Riemannian metrics
via a space-time consideration. Some new monotonicity formulae are
derived. As applications of the new local monotonicity formulae,
some local regularity theorems concerning Ricci flow are proved.
\end{abstract}

\maketitle

\section{Introduction}

The mean value theorem for harmonic functions plays an central role
in the theory of harmonic functions. In this article we discuss its
generalization on manifolds and show how such generalizations lead
to various monotonicity formulae. The main focuses of this article
are  the corresponding results for the parabolic equations,  on
which there have been many works, including  \cite{Fu, W, FG, GL1,
E},  and the application of the new monotonicity formula to the
study of Ricci flow.

Let us start with the Watson's mean value formula \cite{W} for the
heat equation. Let $U$ be a open subset of $\mathbb{R}^n$ (or a
Riemannian manifold). Assume that $u(x, t)$ is a $C^2$ solution to
the heat equation in a parabolic region $U_T=U\times(0, T)$. For any
$(x, t)$ define the `heat ball' by
$$
E(x, t; r):=\left\{ (y, s)\,| s\le t,
\frac{e^{-\frac{|x-y|^2}{4(t-s)}}}{(4\pi (t-s))^{\frac{n}{2}}}\ge
r^{-n}\right\}.
$$
Then
$$
u(x, t)=\frac{1}{r^n}\int_{E(x, t; r)}u(y,
s)\frac{|x-y|^2}{4(t-s)^2}\, dy\, ds
$$
for each $E(x,t; r)\subset U_T$. This result plays an important role
in establishing the Wiener's criterion for heat equation of
$\mathbb{R}^n$ \cite{EG1}.

In \cite{FG}, the above result was generalized to linear parabolic
equations of divergence form. In fact a mean value theorem in terms
of `heat spheres' $\partial E(x, t; r)$ was also derived. This was
later also applied in establishing  the Wiener's criterion for
linear parabolic equations of divergence form on $\mathbb{R}^n$
\cite{GL1}.

More recently, Ecker \cite{E} made the remarkable discovery that
similar mean value property (in terms of  `heat balls') can be
established for the heat equation coupled with the (nonlinear) mean
curvature flow equation. In particular, he discovered a new local
monotonicity formulae for the mean curvature flow through this
consideration.

Since Perelman's celebrated paper \cite{P}, it becomes more evident
that new monotonicity formulae also play  very important roles in
the study of Ricci flow (various monotonicity quantities, such as
isoperimetric constant, entropy like quantity, etc, were  used in
Hamilton's work before). In a recent joint work \cite{EKNT}, the
authors proved a mean value theorem (in terms of `heat balls') for
the heat equation coupled with rather arbitrary deformation equation
on the Riemannian metrics. This in particular gives a local
monotonicity formula for Ricci flow, thanks to the new `reduced
distance' $\ell$-function discovered by Perelman.

In view of the work of \cite{FG}, it is natural to look for mean
value theorem in terms of `heat spheres'  under this general
setting. This is one of the results of this paper. The new mean
value theorem in terms of `heat spheres' also leads new local
monotonicity formulae. The new mean value theorem and monotonicity
formula (in terms of `heat spheres') imply the previous results in
terms of `heat balls' via the integration, as expected. The proof of
the result is through a space-time consideration, namely via the
study of the geometry of $\widetilde{U}=U\times (0, T)$  with the
metric $\tilde{g}(x, t)=g_{ij}(x, t) dx^i dx^j +dt^2$, where $g(x,
t)$ are evolving (by some parabolic equation, such as Ricci flow)
metrics on $\widetilde{U}$. This space-time geometry is considered
un-natural for the parabolic equations since it is not compatible
with parabolic scaling. More involved space-time consideration was
taken earlier by Chow and Chu in \cite{CC}. Nevertheless, it worked
well here for our purpose of proving the mean value theorem and the
monotonicity formulae.

Another purpose of this paper is to derive general mean value
theorem for harmonic functions on Riemannian manifolds. Along this
line there exist several known results before, including those of
Schoen-Yau \cite{SY} for manifolds with nonnegative Ricci curvature,
Michael-Simon \cite{MS} for minimal sub-manifolds, as well as one
for the Cartan-Hadamard manifolds \cite{GW}. Our result here unifies
them all, despite of the simplicity of its derivation. This also
serves the relative easier model for the more involved parabolic
case. Hence we treat it in Section 2 before the parabolic case in
Section 3.

As an application of the new monotonicity formula of \cite{EKNT}, we
formulate and prove a local regularity result for the Ricci flow.
The result and its proof, which are presented in Section 4,  are
motivated by \cite{E}, \cite{Wh} and \cite{P}. This result gives
pointwise curvature estimates under assumption on integral
quantities and hopefully it can shed some lights on suitable
formulation of weak solutions for Ricci flow. One can  refer to
\cite{Ye2} for  results in similar spirits.

It is interesting (also somewhat mysterious) that the quantity used
by P. Li and S.-T. Yau in their fundamental paper \cite{LY} on the
gradient estimates of positive solutions to the heat equation, also
appears in the local monotonicity formulae derived in this article.
Recall that Li-Yau proved that if $M$ is a complete Riemannian
manifold with non-negative Ricci curvature, for any positive
solution $u(x, \tau)$ to the heat equation $\frac{\partial}{\partial
\tau}-\Delta$,
$$
\mathcal{Q}(u):=|\nabla \log u|^2 -(\log u)_\tau \le
\frac{n}{2\tau}.
$$
Here we show  two monotonicity formulae for Ricci flow, which all
involve the expression $\mathcal{Q}$. To state the result we need to
introduce the notion of  the `reduced distance' $\ell$-function of
Perelman. Let $\ell^{(x_0, t_0)}(y, \tau)$ be the `reduced distance'
centered at $(x_0, t_0)$,  with $\tau=t_0-t$. (See \cite{P}, also
Section 4, for a definition.) Define the `sub-heat kernel'
$\hat{K}(y, \tau; x_0, t_0)=\frac{e^{-\ell(y, \tau)}}{(4\pi
\tau)^{\frac{n}{2}}}$. Let $\hat{E}(x_0, t_0; r)$ be the `heat
balls' defined  in terms of  $\hat{K}$ instead. Also define the
`reduced volume' $\theta(\tau):=\int_M \hat{K}\, d\mu(y)$. Then we
have the following result.

\begin{theorem}
Let $(\widetilde{U}, g(t))$ be a solution to Ricc flow. Define $$
\hat{J}(r):=\frac{1}{r^n}\int_{\partial \hat{E}(x_0, t_0;
r)}\frac{\mathcal{Q}(\hat {K})}{\sqrt{|\nabla \log \hat
{K}|^2+|(\log \hat{K})_{\tau}|^2}}\, d\tilde A\,
$$
and
$$
\hat{I}(a, r):=\frac{1}{r^n-a^n}\int_{\hat{E}(x_0, t_0; r)\setminus
\hat{E}(x_0, t_0; a) }\mathcal{Q}(\hat{K})\, d\mu\, d\tau
$$
for $\hat{E}(x_0, t_0; r)\subset \widetilde U$. Here $d\tilde A$ is
the area element of $\partial \hat{E}(x_0, t_; r)$ with respect to
$\tilde g$. Then $\hat{J}(r)\le \hat{I}(a, r)$,  both $\hat{J}(r)$
and $\hat{I}(a, r)$ are monotone non-increasing in $r$. $\hat{I}(a,
r)$ is also monotone non-increasing in $a$. If
$\hat{J}(r_2)=\hat{J}(r_1)$ for some $r_2>r_1$, then $g(t)$ is a
gradient shrinking soliton in $\hat{E}(x_0, t_0; r_2)\setminus
\hat{E}(x_0, t_0; r_1)$ (in fact it satisfies that
$R_{ij}+\ell_{ij}-\frac{1}{2\tau}g_{ij}=0$). If $\hat{I}(a,
r_2)=\hat{I}(a, r_1)$ then $g(t)$ is a gradient shrinking soliton on
$\hat{E}(x_0, t_0; r_2)\setminus \hat{E}(x_0, t_0; a)$. Moreover on
a gradient shrinking soliton, if $t_0$ is the terminating time, then
both $\hat{I}(a, r)$ and $\hat{J}(r)$ are constant and equal to the
`reduced volume' $\theta(\tau)$.
\end{theorem}

The quantity $\mathcal{Q}(\hat K)$ can also be expressed in terms of
the trace LYH expression (modeling the gradient shrinking solitons).
Please see (\ref{ly-rcf}) for details.    The part on $\hat{I}(a,
r)$ of the above result, in the special case $a=0$, has been
established earlier in \cite{EKNT}. (Here our discussion focuses on
the smooth situation while \cite{EKNT} allows  the Lipschitz
functions.) The similar result holds for the mean curvature flow and
the monotonicity of $J$-quantity in terms of `heat spheres' gives a
new local monotone (non-decreasing) quantity for the mean curvature
flow. Please see Section 3 and 5 for details. There is a still open
question of ruling out the grim reaper (cf. \cite{wh2} for an
affirmative answer in the mean-convex case) as a possible
singularity model for the finite singularity of mean curvature flow
of embedded hypersurfaces. It is interesting to find out whether or
not the new monotone quantities of this paper can play any role in
understanding this problem.

\medskip

{\it Acknowledgement}.  We would like to thank Reiner Sch\"atzle,
Jiaping Wang for helpful discussions on two separate technical
issues, Ben Chow for bringing our attention to Schoen-Yau's mean
value inequality, Nicola Garofalo, Richard Hamilton for their
interests and encouragement. We are also grateful to Klaus Ecker for
insisting that we publish the main result of section 4 alone.

\section{The mean value theorem for Laplace equation }
In this section we shall derive a mean value theorem for harmonic
functions  on a fairly large class of Riemannian manifolds. Let $(M,
g)$ be a Riemannian manifold. Let $\Omega$ be a bounded domain in
it. Denote by $G_{\Omega}(x,y)$ the Green's function with Dirichlet
boundary condition on $\partial \Omega$. By the definition
$$
\D_y G_{\Omega}(x, y)=-\delta_x(y).
$$
By the maximum principle $G_{\Omega}(x, y)>0$ for any $x, y \in
\Omega$. By Sard's theorem we know that for almost every $r$, the
`$G_\Omega$-sphere'
$$\Psi_r:=\{y\, | G_\Omega(x,y)= r^{-n}\}$$
is a smooth hypersurface in $\Omega$. Let $\phi_r(y)=G_{\Omega}(x,
y)-r^{-n}$ and let $$ \Omega_r=\{ y\, | \phi_r(y)>0\}
$$
be the `$G_\Omega$-ball'.  We have the following result.

\begin{proposition}\label{mv-general} Let $v$ be a smooth function
on $\Omega$. For every $r>0$
\begin{equation}\label{mv-2}
v(x)=\frac{1}{r^n}\int_{\Omega_r} |\nabla \log G_\Omega|^2 v\,
d\mu-\frac{n}{r^n}\int_0^r \eta^{n}\int_{\Omega_\eta} \phi_\eta
\Delta v\, d\mu \frac{d\eta}{\eta}.
\end{equation}
\end{proposition}
\begin{proof}
We first show that  for almost every $r>0$
\begin{equation}\label{mv-1}
v(x)= \int_{\Psi_r}|\nabla  G_{\Omega}| v\, dA_y
-\int_{\Omega_r}\phi_r \Delta v\, d\mu_y.
\end{equation}
By the Green's second identity
$$
\int_{\Omega_r}\left((\Delta G_\Omega)v -G_\Omega(\Delta v)\right)\,
d\mu =\int_{\Psi_r} \left(\frac{\p G_\Omega}{\p \nu} v-\frac{\p
v}{\p \nu}G_\Omega\right)\, dA
$$
we have that
\begin{equation}\label{11}
v(x)=-\int_{\Psi_r} \left(\frac{\p G_\Omega}{\p \nu} v-\frac{\p
v}{\p \nu}G_\Omega\right)\, dA-\int_{\Omega_r}G_\Omega(\Delta v)\,
d\mu.
\end{equation}
Notice that on $\Psi_r$
\begin{equation}\label{12}
\frac{\p G_{\Omega}}{\p \nu}=-|\nabla G_\Omega|
\end{equation}
and
\begin{eqnarray}\label{13}
\int_{\Psi_r}\frac{\p v}{\p \nu} G_\Omega\,
dA&=&\frac{1}{r^n}\int_{\Psi_r}\frac{\p v}{\p \nu} \, dA\\
&=& \frac{1}{r^n}\int_{\Omega_r}\Delta v\, d\mu.\nonumber
\end{eqnarray}
The equation (\ref{mv-1}) follows by combing (\ref{11}), (\ref{12})
and (\ref{13}).

The equality (\ref{mv-2}) follows from (\ref{mv-1}) by the co-area
formula. In deed, multiplying $\eta^{n-1}$ on the both sides of
(\ref{mv-1}),  integrate on $[0, r]$. Then we have that
\begin{eqnarray*}
\frac{1}{n}r^n v(x)&=&\int_0^r\eta^{n-1}\int_{G_\Omega
=\eta^{-n}}|\nabla G_\Omega|v\, dA\,d\eta+ \int_0^r \eta^n
\int_{\Omega_\eta}\phi_\eta \Delta v\, d\mu \frac{d\eta}{\eta}\\
&=& \frac{1}{n}\int_{r^{-n}}^\infty \frac{1}{\alpha^2}\int_{G_\Omega
=\a}|\nabla G_\Omega|v\, dA\, d\a +\int_0^r \eta^n
\int_{\Omega_\eta}\phi_\eta \Delta v\, d\mu \frac{d\eta}{\eta}\\
&=&\frac{1}{n}\int_{r^{-n}}^\infty\int_{G_\Omega =\a} \frac{|\nabla
\log G_\Omega|^2 v}{|\nabla G_\Omega|}\, dA \, d\a+\int_0^r \eta^n
\int_{\Omega_\eta}\phi_\eta \Delta v\, d\mu \frac{d\eta}{\eta}\\
&=&\frac{1}{n}\int_{\Omega_r} |\nabla \log G_\Omega|^2 v\, d\mu
+\int_0^r \eta^{n}\int_{\Omega_\eta} \phi_\eta \Delta v\, d\mu
\frac{d\eta}{\eta}.
\end{eqnarray*}
In the last equation we have applied the co-area formula \cite{EG},
for which we have to verify that $|\nabla \log G_\Omega|^2$ is
integrable. This follows from the asymptotic behavior of
$G_\Omega(x, y)$ near $x$. It can also be seen, in the case of $n\ge
3$, from a known estimate (cf. \cite{Li} Theorem 6.1) of Cheng-Yau,
which  asserts that near $x$, $ |\nabla \log G_\Omega|^2(y)\le
A\left(1+\frac{1}{d^2(x, y)}\right). $ for some $A>0$.
\end{proof}

It would be more convenient to apply the Proposition
\ref{mv-general} if we can replace the Green's function $G_\Omega(x,
y)$, which depends on $\Omega$ and usually hard to estimate/compute,
by a canonical positive Green's function.
 To achieve this we need to assume that the
Riemannian manifold $(M, g)$ is non-parabolic. Namely, there exists
a minimum positive Green's function $G(x, y)$ on $M$. However, on
general non-parabolic manifolds, it may happen that $G(x, y)$ does
not approach to $0$ as $y\to \infty$. This would imply that the
integrals on the right hand side of (\ref{mv-1}) or (\ref{mv-2}) are
over a noncompact domain (a unbounded hyper-surface). Certain
requirements on $v$ (justifications) are needed to make sense of the
integral. It certainly works for the case that $v$ has compact
support. If we want to confine ourself with the situation that
$\Omega_r=\{y\, |G(x, y)\ge r^{-n}\}$ is compact we have to impose
further that
\begin{equation}\label{snp} \lim_{y\to \infty}G(x, y)=0.
\end{equation}

\begin{definition}\label{SNP} We call the Riemannian manifold $(M,
g)$ is \textit{ strongly non-parabolic} if (\ref{snp}) holds for the
minimal positive Green's function.
\end{definition}

Now similarly we can define the `$G$-sphere' $\Psi_r$, the
`$G$-ball' $\Omega_r$, and the function $\phi_r$ as before. By
verbatim repeating the proof of Proposition \ref{mv-general}   we
have the following global result.

\begin{theorem}\label{mv-snp} Assume that $(M, g)$ is a strongly
non-parabolic Riemannian manifold. Let $v$ be a smooth function.
Then for every $r>0$
\begin{equation}\label{mv-4}
v(x)=\frac{1}{r^n}\int_{\Omega_r} |\nabla \log G|^2 v\, d\mu
-\int_0^r \frac{n}{\eta^{n+1}}\int_{\Omega_\eta} \psi_\eta \Delta
v\, d\mu d\eta.
\end{equation}
Here $\psi_r=\log(Gr^n)$. Also for almost every $r>0$
\begin{equation}\label{mv-3}
v(x)= \int_{\Psi_r}|\nabla  G| v\, dA_y -\int_{\Omega_r}\phi_r
\Delta v\, d\mu_y.
\end{equation}
\end{theorem}
\begin{proof}
The only thing we need to verify is that
$$
\int_0^r \frac{n}{\eta^{n+1}}\int_{\Omega_\eta} \psi_\eta \Delta v\,
d\mu d\eta=\frac{n}{r^n}\int_0^r
\eta^{n}\int_{\Omega_\eta}\phi_\eta\Delta v
 \, d\mu \frac{d\eta}{\eta}
$$
which can be checked directly. Please see also (\ref{equ1}) for a
more  general equality.
\end{proof}

\begin{remark}\label{rk1}
a) If $(M, g)=\R^n$ with $n\ge 3$, then  $$G(x,
y)=\frac{1}{n(n-2)\omega_n}d^{2-n}(x, y),$$ with $\omega_n$ being
the volume of unit Euclidean ball and $d(x, y)$ being the distance
between $x$ and $y$. Clearly $M$ is strongly non-parabolic. If $v$
is a harmonic function, a routine exercise shows that  Theorem
\ref{mv-snp} implies the classical mean value theorem
$$
v(x)= \frac{1}{nR^{n-1}\omega_n}\int_{\partial B(x, R)}v(y)\, dA.
$$

b) The regularity on $v$ can be weaken to, say being  Lipschitz. Now
we have to understand $\Delta v$ in the sense of distribution and
the integral $\int_{\Omega_r} \phi_r \Delta v$ via a certain
suitable approximation \cite{E, EKNT, N3}.

c) Li and Yau  \cite{LY} proved that if $(M, g)$ has nonnegative
Ricci curvature, then $(M, g)$ is non-parabolic if and only if
$\int_r^\infty \frac{\tau}{V_x(\tau)}\, d\tau <\infty$, where
$V_x(\tau)$ is the volume of ball $B(x, \tau)$. Moreover, there
exists $C=C(n)$ such that
$$
C^{-1}\int_{d(x, y)}^\infty \frac{\tau}{V_x(\tau)}\, d\tau\le G(x,
y)\le C\int_{d(x, y)}^\infty \frac{\tau}{V_x(\tau)}\, d\tau.
$$
This shows that any non-parabolic $(M, g)$ is  strongly
non-parabolic. In this case by the gradient estimate of Yau
\cite{Li} we have that there exists $C_1(n)$ such that if  $\D v\ge
0$ and $v\ge 0$
\begin{equation}\label{toLiSch}
v(x)\le \frac{C_1(n)}{r^n}\int_{\Psi_r}\frac{v(y)}{d(x, y)}\, dA_y.
\end{equation}
One would expect that this should imply the well-known mean value
theorem of Li-Schoen \cite{LS}.

d) A example class of strongly non-parabolic Riemannian manifolds
 are those manifolds with the  Sobolev inequality:
 \begin{equation}\label{sobolev}
\left(\int_M f^{\frac{2\nu}{\nu -2}}\,
d\mu\right)^{\frac{\nu-2}{\nu}}\le A\int_M |\nabla f|^2\, d\mu
 \end{equation}
for some $\nu>2$, $A>0$, for any smooth $f$ with compact support. In
deed, by \cite{Davies} and \cite{Gri}, the Sobolev inequality
(\ref{sobolev}) implies that there exist $B$ and $D>0$ such that
$$
H(x,y, t)\le B t^{-\frac{\nu}{2}}\exp\left(-\frac{d^2(x, y)}{D
t}\right).
$$
This implies the estimate
$$
G(x, y)\le C(B, \nu, D) d^{-\nu+2}(x, y)
$$
for some $C$. Therefore $(M, g)$ is strongly non-parabolic.

e) Another set of examples are the Riemannian manifolds with
positive lower bound on the spectrum with respect to the Laplace
operator \cite{LW}. On this type of manifolds we usually   assume
that the Ricci curvature $Ric\ge -(n-1) g$. Since the manifold may
contain ends with finite volume (as examples in \cite{LW} show), on
which the Green's function  certainly does not tend to zero at
infinity, further assumptions are needed to ensure the strongly
non-parabolicity. If we assume further that there exists a
$\delta>0$ such that $V_x(1)\ge \delta$, we can show that it is
strongly non-parabolic (shown in the proposition below). Note that
the assumption holds for the universal covers of  compact Riemannian
manifolds with $Ric\ge -(n-1)g$.

\end{remark}

\begin{proposition}\label{ps} Assume that $(M^n, g)$ is a complete
Riemannian manifold with positive lower bound on the spectrum of the
Laplace operator. Suppose that $Ric\ge -(n-1) g$ and that there
exists a $\delta>0$ such that $V_x(1)\ge \delta$ for all $x\in M$.
Then $(M, g)$ is strongly non-parabolic.
\end{proposition}
\begin{proof}
First we have the following upper bound of the heat kernel
\begin{equation}\label{Li-lhf}
H(x, y, t)\le C_1\exp(-\lambda
t)V^{-\frac{1}{2}}_x(\sqrt{t})V^{-\frac{1}{2}}_y(\sqrt{t})\exp\left(-\frac{d^2(x,
y)}{Dt}+C_2\sqrt{t}\right)
\end{equation}
for some absolute constant $D>4$, $C_2=C_2(n)$ and $C_1=C_1(D, n)$.
Please see \cite{Li2} for a proof. Here $\lambda>0$ denotes the
greatest lower bound on the spectrum of the Laplace operator. The
assumption on the lower bound of the curvature implies that
$$
V^{-\frac{1}{2}}_x(\sqrt{t})V^{-\frac{1}{2}}_y(\sqrt{t})\le h(t)
$$
where
\[
h(t)=\left\{ \begin{array}{ll} \frac{1}{\delta} & \mbox{if }t\ge1\\
C(\delta, n)t^{-\frac{n}{2}} & \mbox{if } t\le 1 \end{array}.\right.
\]
By some elementary computation and estimates we have that
\begin{eqnarray*}
G(x, y)&=&\int_0^{\infty}H(x, y, t)\, dt\\
&\le & \frac{4}{\delta
\lambda}e^{\frac{C_2^2}{2\lambda}}\exp\left(-\sqrt{\frac{\lambda}{D}}d(x,
y)\right)\\
&\quad &+ C(\delta, n, D) e^{\frac{C_2^2}{2\lambda}}d^{-n+2}(x,
y)\int_{\frac{d^2(x, y)}{D}}^\infty
\exp(-\tau)\tau^{\frac{n}{2}-2}\, d\tau
\end{eqnarray*}
which goes to $0$ as $d(x, y)\to \infty$. In fact one can have the
upper bound $ G(x, y)\le C\exp\left(-\sqrt{\frac{\lambda}{D}} d(x,
y)\right)$.
\end{proof}

There exists  monotonicity formulae related to Proposition
\ref{mv-general} and Theorem \ref{mv-snp}, which we shall illustrate
below. For simplicity we just denote by $G$ for both $G_\Omega(x,
y)$, in the  discussion concerning a bounded domain $\Omega$, and
$G(x, y)$, when it is on the strongly non-parabolic manifolds.

For any smooth (or Lipschitz) $v$ we define
$$
I_v(r):=\frac{1}{r^n}\int_{\Omega_r} |\nabla \log G|^2 v\, d\mu
$$
and
$$
J_v(r):=\int_{\Psi_r}|\nabla G|v\, dA.
$$
They are related through the relation $r^n I_v(r)=n \int_0^r
\eta^{n-1} J(\eta)\, d\eta$, which can be shown by using  Tonelli's
theorem and the co-area formula.

\begin{corollary}\label{mono1}
Let $\psi_r=\log(Gr^n)$, which is nonnegative on $\Omega_r$. We have
that for almost every $r>0$
\begin{equation}\label{mono-e1}
\frac{d }{d r} I_v(r)=\frac{n}{r^{n+1}}\int_{\Omega_r}\left(\D
v\right)\psi_r\, d\mu;
\end{equation}
\begin{equation}\label{mono-s1}
\frac{d}{d r} J_v(r)=\frac{n}{r^{n+1}}\int_{\Omega_r}\D v\, d\mu.
\end{equation} In particular, $I_v(r)$ and $J_v(r)$ are monotone non-increasing
(non-decreasing) in $r$, provided that  $v$ is super-harmonic
(sub-harmonic).
\end{corollary}
\begin{proof}
Differentiate (\ref{mv-4}) with respect to $r$. Then
$$
I'_v(r)=-\frac{n^2}{r^{n+1}}\int_0^r\eta^n\int_{\Omega_\eta}
\phi_\eta \Delta v\, d\mu
\frac{d\eta}{\eta}+\frac{n}{r}\int_{\Omega_r}\phi_r \D v\, d\mu.
$$
On the other hand, Tonelli's theorem gives
\begin{eqnarray*}
\int_0^r\eta^n\int_{\Omega_\eta} \phi_\eta \Delta v\, d\mu
\frac{d\eta}{\eta}&=& \int_{G\ge
r^{-n}}\int_{G^{-\frac{1}{n}}}^r\left( G
\eta^{n-1}-\frac{1}{\eta}\right)\D v\, d\mu \, d\eta\\
&=& \frac{r^n}{n}\int_{\Omega_r} (\D v)(G-r^{-n})\, d\mu
-\frac{1}{n}\int_{\Omega_r}(\Delta v) \psi_r\, d\mu.
\end{eqnarray*}
From the above two equations, we have the claimed (\ref{mono-e1}),
observing that $\phi_r=G-r^{-n}$. The proof of (\ref{mono-s1}) is
very similar.
\end{proof}

A by-product of the above proof is that for any $f(y)$ (regular
enough to makes sense the integrals)
\begin{equation}\label{equ1}
\frac{n}{r^n}\int_0^r \eta^{n}\int_{\Omega_\eta}f \phi_\eta \, d\mu
\frac{d\eta}{\eta}=\int_0^r \frac{n}{\eta^{n+1}}\int_{\Omega_\eta}
f\psi_\eta\, d\mu\, d\eta.
\end{equation}
In fact  let $F(r)$ be the left hand side, the above proof shows
that
$$
F'(r)=\frac{n}{r^{n+1}}\int_{\Omega_r}f\psi_r \, d\mu.
$$

In \cite{SY}, Schoen and Yau proved a mean value inequality for
sup-harmonic functions with respect to sufficiently small  geodesic
balls, when $(M, g)$ is a complete Riemannian manifold (or just a
piece of) with nonnegative Ricci curvature. More precisely, they
showed that if $x\in M$ and $B(x, R)$ lies inside a normal
coordinate centered at $x$. Let $f\ge0$ be a Lipschitz function
satisfying that $\Delta f\le 0$ then
\begin{equation}\label{sy}
f(x)\ge \frac{1}{n\omega_n R^{n-1}}\int_{\partial B(x, R)}f(y)\, dA.
\end{equation}
We shall see that a consequence of Theorem \ref{mv-snp} implies a
global version of the above result.  First recall the following
fact. If $(M, g)$ is a Riemannian manifold such that its Ricci
curvature satisfies $Ric\ge -(n-1)k^2 g$, then
\begin{equation}\label{green-com}
\Delta_y\hat {G}(d(x, y))\ge -\delta_x(y).\end{equation} We call
such $\hat{G}$ a sub-Green's function. Here $\hat{G}(\bar{d}(x, y))$
is the Green's function of the space form $(\overline{M}, \bar{g})$
with constant curvature $-k^2$, where $\bar{d}$ is the distance
function of $\overline M$. One can refer to Chapter 3 of \cite{SY}
for a proof of the corresponding parabolic result (which is
originally proved in \cite{CY}). The key facts used for the proof
are the standard Laplace comparison theorem \cite{SY} and
$\hat{G}'\le 0$ (which follows from the maximum principle). Now we
can define
$$
\hat{\Omega}_r=\{y\, | \hat{G}\ge r^{-n}\}, \quad \quad
\hat{\Psi}_r=\{y\, |\hat{G}=r^{-n}\}
$$
and let $\hat{\phi}_r=\hat{G}-r^{-n}$.

 The virtue of the proof of Proposition
\ref{mv-general} gives the following  general result.

\begin{proposition}\label{mv-rcbb} Let $(M, g)$ be a complete
Riemannian manifold with $Ric\ge -(n-1)k^2g$. In the case $k=0$ we
assume further that $n\ge 3$. Let $v\ge 0$ be a Lipschitz function.
Then
\begin{equation}\label{mv-5}
v(x)\ge  \int_{\hat{\Psi}_r}|\nabla  \hat{G}| v\, dA_y
-\int_{\hat{\Omega}_r}\hat{\phi}_r \Delta v\, d\mu_y.
\end{equation}
For every $r>0$
\begin{equation}\label{mv-6}
v(x)\ge \frac{1}{r^n}\int_{\hat{\Omega}_r} |\nabla \log \hat{G}|^2
v\, d\mu -\frac{n}{r^n}\int_0^r \eta^{n}\int_{\hat{\Omega}_\eta}
\hat{\phi}_\eta \Delta v\, d\mu \frac{d\eta}{\eta}.
\end{equation}
\end{proposition}

By some straight forward computations one can verify  that
(\ref{mv-5}), in the case of $k=0$,  implies a global version of
(\ref{sy}).

Moreover, in this non-exact case there are still monotonicity
formulae related to Proposition \ref{mv-rcbb}.

\begin{corollary}\label{mono2}
Let
$$
\hat{I}_v(r)= \frac{1}{r^n}\int_{\hat{\Omega}_r} |\nabla \log
\hat{G}|^2 v\, d\mu, \quad  \hat{J}_v(r)=\int_{\hat{\Psi}_r} |\nabla
\hat{G}| v\, dA.
$$
 If $v\ge 0$, then for almost every $r$,
\begin{equation}\label{mono-e2}
\frac{d}{d r} \hat{I}_v(r)\le
\frac{n}{r^{n+1}}\int_{\hat{\Omega}_r}(\Delta v)\hat{\psi}_r\, d\mu
\end{equation}
where $\hat{\psi}_r=\log(\hat{G}r^n)$, and
\begin{equation}\label{mono-s2}
\frac{d}{d r}\hat{J}_v(r)\le\frac{n}{r^{n+1}}\int_{\hat{\Omega}_r}\D
v\, d\mu.
\end{equation}
In particular, if $\D v\le 0$, $\hat{I}_v(r)$ and $\hat{J}_v(r)$ are
monotone non-increasing in $r$. Moreover, if the equality holds in
the inequality (\ref{mono-e2}) (or (\ref{mono-s2})) for some $v>0$
at some $r>0$ then $B(x, R)$, the biggest ball contained in
$\hat{\Omega}_r$, is isometric to the corresponding ball in the
space form.
\end{corollary}
\begin{proof} Here we adapt a different scheme from that of
Corollary \ref{mono1}. For the convenience we let $\hat{\psi}=\log
\hat{G} $. Differentiate $\hat I_v(r)$ we have
\begin{equation}\label{deri}
\hat{I}_v'(r)=-\frac{n}{r^{n+1}}\int_{\hat{\Omega}_r}|\nabla
\hat{\psi}_r|^2 v\, d\mu
+\frac{n}{r^{n+1}}\int_{\hat{\Psi}_r}|\nabla \hat{\psi}|v\, dA.
\end{equation}
Using that $\frac{\p}{\p \nu} \hat{\psi}=-|\nabla \hat{\psi}|$ on
$\hat{\Psi}_r$, we have that
\begin{eqnarray*}
\int_{\hat{\Psi}_r}|\nabla \hat{\psi}|v\,
dA&=&-\int_{\hat{\Psi}_r}\left(\frac{\p}{\p \nu} \hat{\psi}\right) v\, dA\\
&=&\int_{\hat{\Omega}_r} \left((\D v)\hat{\psi}-(\Delta
\hat{\psi})v\right)\, d\mu -\int_{\hat{\Psi}_r}\left(\frac{\p}{\p
\nu} v\right) \hat{\psi}\, dA\\
&=&\int_{\hat{\Omega}_r} \left((\D v)\hat{\psi}-(\Delta
\hat{\psi})v\right)\, d\mu+\log r^n \int_{\hat{\Omega}_r}\D v\,
d\mu.
\end{eqnarray*}
Using (\ref{green-com}), and the fact that $\frac{1}{\hat {G}}\to 0
$ as $y\to x$, we have that
$$
-\int_{\hat{\Omega}_r}(\Delta \hat{\psi})v\, d\mu\le
\int_{\hat{\Omega}_r}|\nabla \hat{\psi}|^2v\, d\mu.
$$
Therefore
$$ \int_{\hat{\Psi}_r}|\nabla \hat{\psi}|v\, dA\le\int_{\hat{\Omega}_r}|\nabla \hat{\psi}|^2v\,
d\mu+\int_{\hat{\Omega}_r}(\Delta v)\hat{\psi}_r\, d\mu.
$$
Together with (\ref{deri}) we have the claimed inequality
(\ref{mono-e2}). The equality case follows from the equality in the
Laplace comparison theorem.

To prove (\ref{mono-s2}), for any $r_2>r_1$, the Green's second
identity implies that
\begin{eqnarray*}
J_v(r_2)-J_v(r_1)&=&\int_{\hat{\Psi}_{r_2}}\left(-\frac{\p
\hat{G}}{\p \nu}\right)v\, dA +\int_{\hat{\Psi}{r_1}}\frac{\p
\hat{G}}{\p \nu}
v\, dA\\
&=& -\int_{\hat{\Omega}_{r_2}\setminus \hat{\Omega}_{r_1}}\left((\D
\hat{G})v -(\D v) \hat{G}\right)\, d
\mu\\
&\quad& -\int_{\hat{\Psi}_{r_2}}G\frac{\p v}{\p \nu}\, dA
+\int_{\hat{\Psi}_{r_1}}G\frac{\p v}{\p \nu}\, dA\\
&\le & \int_{\hat{\Omega}_{r_2}}\D
v(\hat{G}-r_2^{-n})-\int_{\hat{\Omega}_{r_1}}\D v
(\hat{G}-r_1^{-n})\, d\mu.
\end{eqnarray*}
Dividing by $(r_2-r_1)$, the claimed result follows by co-area
formula and taking limit $(r_2-r_1)\to 0$.
\end{proof}

The general result can also be applied to a somewhat opposite
situation, namely to Cartan-Hardamard manifolds and  minimal
sub-manifolds in such manifolds. Recall that $(M^n, g)$  is called a
Cartan-Hardamard manifold if it is simply-connected with the
sectional curvature $K_M\le 0$. Let $ \tilde G(x,
y)=\frac{1}{n(n-2)\omega_n}d^{2-n}(x, y)$ ($n\ge 3$). Then we have
that \cite{GW}
\begin{equation}\label{dual}
\D_y  \tilde{G}(x, y)\le -\delta_x(y).
\end{equation}
Similarly we call such $\tilde G$ a sup-Green's function. Let $N^k$
($k\ge 3$) be a minimal (immersed) submanifold in $M^n$ (where $M$
is the Cartan-Hardamard manifold as the above. Let $\bar{\Delta}$ be
the Laplace operator of $N^k$. It is easy to check that
\begin{equation}\label{dual2}
\bar{\D} \bar{G}\le -\delta_x(y)
\end{equation}
where $\bar{G}(x, y)=\frac{1}{n(n-2)\omega_k}d^{2-k}(x, y)$. Here
$d(x, y)$ is the extrinsic distance function of $M$. Similarly we
can define $\tilde{\Omega}_r\subset M$ for the first case and
$\overline{\Omega}_r\subset N$ for the minimum submanifold case.

\begin{proposition}\label{mv-ch} Let $(M, g)$ be Cartan-Hardamard manifold.  Let $v\ge 0$ be a Lipschitz function.
Then
\begin{equation}\label{mv-51}
v(x)\le \int_{\tilde{\Psi}_r}|\nabla  \tilde{G}| v\, dA_y
-\int_{\tilde{\Omega}_r}\tilde{\phi}_r \Delta v\, d\mu_y.
\end{equation}
For every $r>0$
\begin{equation}\label{mv-61}
v(x)\le \frac{1}{r^n}\int_{\tilde{\Omega}_r} |\nabla \log
\tilde{G}|^2 v\, d\mu -\frac{n}{r^n}\int_0^r
\eta^{n}\int_{\tilde{\Omega}_\eta} \tilde{\phi}_\eta \Delta v\, d\mu
\frac{d\eta}{\eta}.
\end{equation}
\end{proposition}

We state the result only for $\tilde{G}$ since the exactly same
result holds $\bar{G}$. The mean value inequality (\ref{mv-51}) of
$\bar{G}$ recovers the well-known result of Michael and Simon
\cite{MS} if $N$ is a minimal sub-manifold in the Euclidean spaces.
 See also \cite{CLY} for another proof of Michael-Simon's result via
 a  heat kernel comparison theorem.

\begin{corollary}\label{mono3}
Let
$$
\tilde{I}_v(r)= \frac{1}{r^n}\int_{\tilde{\Omega}_r} |\nabla \log
\tilde{G}|^2 v\, d\mu, \quad \tilde{J}_v(r)=\int_{\tilde{\Psi}_r}
|\nabla \tilde{G}| v\, dA.
$$
 If $v\ge 0$, then for almost every $r$,
\begin{equation}\label{mono-e3}
\frac{d}{d r} \tilde{I}_v(r)\ge
\frac{n}{r^{n+1}}\int_{\tilde{\Omega}_r}(\Delta v)\tilde{\psi}_r\,
d\mu
\end{equation}
where $\tilde{\psi}_r=\log(\tilde{G}r^n)$, and
\begin{equation}\label{mono-s3}
\frac{d}{d
r}\tilde{J}_v(r)\ge\frac{n}{r^{n+1}}\int_{\tilde{\Omega}_r}\D v\,
d\mu.
\end{equation}
In particular, if $\D v\ge 0$, $\tilde{I}_v(r)$ and $\tilde{J}_v(r)$
are monotone non-increasing in $r$. Moreover, if the equality holds
in the inequality (\ref{mono-e3}) (or (\ref{mono-s3})) for some
$v>0$ at some $r>0$ then $B(x, R)$, the biggest ball contained in
$\tilde{\Omega}_r$, is isometric to the corresponding ball in the
Euclidean space.
\end{corollary}

\section{Spherical mean value theorem for heat equation with changing metrics}

In \cite{EKNT},  mean value theorems in terms of `heat balls'  were
proved for heat equations with respect to evolving metrics. A mean
value theorem in terms of `heat spheres' was also established for
the heat equation with respect to a fixed Riemannian metric. Since
the mean value theorems in terms of `heat balls' are usually
consequences  of the ones in terms of `heat spheres' it is desirable
to have the mean value theorem in terms of `heat spheres' for the
heat equation with respect to  evolving metrics too. Deriving such a
mean value theorem is the main result of this section. The proof is
an adaption  of the argument of  \cite{FG} for operators of
divergence form on Euclidean case (see also \cite{Fu}) to the
evolving metrics case. The key is the Green's second identity
applied to the space-time.

 Let $(M,
g(t))$ be a family of metrics evolved by the equation
\begin{equation}\label{evolution}
\frac{\p}{\p t} g_{ij}=-2\Upsilon_{ij}.\end{equation}We start by
setting up some basic space-time notions. First we fix a point
$(x_0, t_0)$ in the space time. For the simplicity we assume that
$t_0=0$. Assume that $g(t)$ is a solution to (\ref{evolution}) on
$M\times (\alpha, \beta)$ with $\a<0<\beta$. Denote $M\times (\a,
\b)$ by $\widetilde{M}$, over which we define  the metric $\tilde
g(x,t)=g_{ij}(x, t)dx^i dx^j + dt^2$, where $t$ is the global
coordinate of $(\a, \b)$. We consider the heat operator $\heat$ with
respect to the changing metric $g(t)$. Now the {\it conjugate heat
operator} is $ \frac{\partial }{\p t}+\D -R(y, t)$, where $R(y,
t)=g^{ij}\Upsilon_{ij}$. We need some elementary space-time
computations for $(\widetilde{M}, \tilde g)$. In \cite{CC} (see also
\cite{P}), a similar consideration was originated, but for some
degenerate metrics satisfying the generalized Ricci flow equation
instead. In the following  indices $i,j, k$ are between $1$ and $n$,
$A, B, C$ are between $0$ and $n$. The index $0$ denotes the $t$
direction.

\begin{lemma}\label{connection}
Let $\Gamma_{ij}^k$ be the Christoffel symbols of $g_{ij}$, and let
$\widetilde{\Gamma}_{AB}^C$ be  the ones for $\widetilde g$. Then
\begin{eqnarray}\label{conn}
&\, &\widetilde{\Gamma}^{0}_{ij}=\Upsilon_{ij},\\
&\, &\widetilde{\Gamma}^{0}_{0, k}=\widetilde{\Gamma}^{0}_{00}=0,\\
&\, &\widetilde{\Gamma}^{k}_{ij}=\Gamma_{ij}^k,\\
&\,&\widetilde{\Gamma}^i_{0k}=-\widetilde{g}^{il}\Upsilon_{lk}=-\Upsilon^{i}_k.
\end{eqnarray}
\end{lemma}
\begin{proof}
It  follows from  straight forward computations.
\end{proof}

\begin{lemma}\label{divergence} Let $X$ be a time dependent vector field on $\widetilde{M}$ (namely for each $t$, $X(t)$ is a vector field of $M$). Let
$\widetilde{X}=X+X^0\frac{\p}{\p t}$ for some function $X^0$ on
$\widetilde{M}$. Let $\widetilde{\operatorname{div}}$ be the
divergence operator with respect to $\tilde g$, and let
$\operatorname{div}$ be the divergence operator with respect to $g$
(on $M\times\{t\}$). Then
\begin{equation}\label{div}
\widetilde{\operatorname{div}}
\widetilde{X}=\operatorname{div}(X)-X^0R+\frac{\p}{\p t} X^0.
\end{equation}
\end{lemma}
\begin{proof} Follows from Lemma \ref{connection} and routine
computations.
\end{proof}

Let $H(y, t; x_0, 0)$ be a fundamental solution to the {\it
conjugate heat operator} centered at $(x_0, 0)$. Let $\tau=-t$. By
the abuse of the notation we sometime also write the fundamental
solution as $H(y, \tau; x_0, 0)$, by which we mean the fundamental
solution  with respect to the operator $\frac{\p}{\p \tau}-\D +R$.
Now we define the `heat ball' by $E_r=\{(y, \tau)\, |\, H( y, \tau;
x_0, 0)\ge r^{-n}\}$. By the asymptotics of the fundamental solution
\cite{GL} we know that $E_r$ is compact for $r$ sufficiently small.
Following \cite{FG} we define $E_r^{s}=\{(y, \tau)\in E_r,\, t<s\}$
and two portions of its boundary $P_1^s=\{(y, \tau)\, |\, H(y, \tau;
x_0, 0)=r^{-n},\, t<s\}$ and $P_2^s=\{(y, \tau)\in
\overline{E^s_r},\, t=s\}$. $P_1^0=\p E_r$ is the `heat sphere',
which is the boundary of the `heat ball'.  Let $\psi_r=H(y, \tau;
x_0, 0)-r^{-n}$. We then have the following spherical mean value
theorem.

\begin{theorem}\label{mv-sph} Let $v$ be a smooth function on $\widetilde{M}$.  Then
\begin{eqnarray}\label{mv-sph-gen}
v(x_0, 0)&=&\int_{\p E_r} v\frac{|\nabla H|^2}{\sqrt{|\nabla
H|^2+|H_t|^2}}\, d\tilde A+\frac{1}{r^n}\int_{E_r}R v\, d\mu\,
dt\nonumber\\
&\quad &+\int_{E_r}\phi_r\left(\frac{\p}{\p t}-\D\right)v  \, d\mu\,
dt.
\end{eqnarray}
Here $d\tilde A$ is the $n$-dimensional measure induced from $\tilde
g$.
\end{theorem}
\begin{proof} This follows as in, for
instance, Fabes-Garofalo \cite{FG}, by applying a divergence theorem
in $\widetilde{M}$ and  the above Lemma \ref{div} for the changing
metrics. More precisely, let $\widetilde{X}= v\nabla H-H\nabla v+Hv
\frac{\p}{\p t}$. Applying the divergence theorem to $\tilde{X}$ on
$E^s_r$ we have that
\begin{eqnarray}\label{div-app1}
\int_{E^s_r}\left(v(\frac{\p}{\p t}+\Delta -R)H+H(\frac{\p}{\p
t}-\D)v\right)\, d\mu\, dt&=&
\int_{E^s_r}\widetilde{\mbox{div}}(\widetilde{X})\, d\mu\, dt \nonumber\\
&=&\int_{\p E^s_r} \langle \widetilde{X}, \tilde{\nu}\rangle d\tilde
A.
\end{eqnarray}
Here $\tilde{\nu}$ is the normal $\p E^s_r$  with respect to
$(\widetilde{M}, \tilde g)$. On $P_1^s$ it is given by
$$
\tilde{\nu} =\left(-\frac{\nabla H}{\sqrt{|\nabla H|^2+|H_t|^2}},
-\frac{H_t}{\sqrt{|\nabla H|^2+|H_t|^2}}\right).
$$
On $P_2^s$ it is just $\frac{\p}{\p t}$. Compute
\begin{eqnarray}
\label{div-app2}\int_{\p E^s_r} \langle \widetilde{X},
\tilde{\nu}\rangle d\tilde A&=&\int_{P_1^s} \left(-v\langle \nabla
H, \tilde{\nu}\rangle -H \langle \nabla v, \tilde{\nu}\rangle + Hv
\langle \frac{\p}{\p t}, \tilde{\nu}\rangle\right)\, d \tilde A\nonumber\\
&\quad&+\int_{P_2^s} vH\, d\mu_s\nonumber\\
&=&\int_{P_1^s} \left(-v\frac{|\nabla H|^2}{\sqrt{|\nabla
H|^2+|H_t|^2}} -H \langle \nabla v,
\tilde{\nu}\rangle\right.\nonumber\\
&\quad & +\left. Hv \langle \frac{\p}{\p t}, \tilde{\nu}\rangle
\right)\, d \tilde A+\int_{P_2^s} vH\, d\mu_s\nonumber\\
&=&\int_{P_1^s} \left(-v\frac{|\nabla H|^2}{\sqrt{|\nabla
H|^2+|H_t|^2}} \right)\, d \tilde A+\int_{P_2^s} vH\, d\mu_s \\
&\quad & +r^{-n}\int_{P_1^s} \left(-\langle \nabla v,
\tilde{\nu}\rangle+v \langle \frac{\p}{\p t},
\tilde{\nu}\rangle\right)\,d\tilde{A}\nonumber
\end{eqnarray}
Similarly, applying the divergence theorem to $\widetilde{Y}=-\nabla
v+v \frac{\p}{\p t}$  gives that
\begin{eqnarray}\label{div-app3}
\int_{E^s_r}\left(-\D v-Rv+\frac{\p v}{\p t}\right)\, d\mu\, dt&=&
\int_{P_1^s} \left(-\langle \nabla v, \tilde{\nu}\rangle+v
\langle \frac{\p}{\p t}, \tilde{\nu}\rangle\right)\, d\tilde{A}\nonumber\\
&\quad & +\int_{P^s_2}v\, d\mu_s.
\end{eqnarray}
By (\ref{div-app1})--(\ref{div-app3}) we have that
\begin{eqnarray*}
\int_{P^s_2}v (H-r^{-n})\, d\mu_s &=& \int_{E^s_r}
\phi_r\left(\frac{\p v}{\p t}-\D v\right)\, d\mu\, dt
+r^{-n}\int_{E^s_r}R v\, d\mu\, dt \\
&\quad &+\int_{P_1^s} \left(v\frac{|\nabla H|^2}{\sqrt{|\nabla
H|^2+|H_t|^2}} \right)\, d \tilde A.
\end{eqnarray*}
Letting $s\to 0$, the claimed result follows from the asymptotics
\begin{equation}\label{claim}
\lim_{s\to 0}\int_{P^s_2}v (H-r^{-n})\, d\mu_s=v(x_0, 0)
\end{equation}
which can be checked directly using the asymptotics of $H$, or
simply the definition of $H$.
\end{proof}

\begin{remark}\label{rk2}
Theorem \ref{mv-sph} is the parabolic analogue of (\ref{mv-1}) and
(\ref{mv-3}). In the spacial case of  the metrics being  fixed,
namely $\Upsilon =0$, Theorem \ref{mv-sph} gives a manifold version
of the spherical mean value theorem for solutions to the heat
equation. Please see   \cite{Fu} and \cite{FG} for earlier results
when $M=\mathbb{R}^n$.

As a corollary of Theorem \ref{mv-sph} one can obtain the `heat
ball' mean value theorem proved in \cite{EKNT} by integrating
(\ref{mv-sph-gen}) for $r$ and applying the co-area formula and
Tonelli's theorem. Indeed multiplying $\eta^{n-1}$ on both sides of
(\ref{mv-sph-gen}), then integrating from $0$ to $r$ as in
Proposition \ref{mv-general},  we have that
\begin{eqnarray}\label{mv-sph-ball}
v(x_0, 0)&=& \frac{1}{r^n}\int_{E_r}\left(|\nabla \log
H|^2+R\psi_r\right)v\, d\mu dt \\
&\quad &+\frac{n}{r^n}\int_0^r \eta^n
\int_{E_\eta}\phi_\eta\left(\frac{\p}{\p t}-\D\right)v  \, d\mu\,
dt\, \frac{d\eta}{\eta}.\nonumber
\end{eqnarray}
By an argument similar to the proof of  (\ref{equ1}) we have that
\begin{eqnarray*}
&\quad &\frac{n}{r^n} \int_0^r \eta^n
\int_{E_\eta}\phi_\eta\left(\frac{\p}{\p t}-\D\right)v  \, d\mu\,
dt\,
\frac{d\eta}{\eta}\\
&=&\int_0^r\frac{n}{\eta^{n+1}}\int_{E_\eta}\psi_\eta\left(\frac{\p}{\p
t}-\D\right)v  \, d\mu\, dt\, d\eta.
\end{eqnarray*}
Therefore (\ref{mv-sph-ball}) is the same as the formula in
\cite{EKNT}. Similar to Corollary \ref{mono1}, we can recover the
monotonicity formulae of \cite{EKNT} from (\ref{mv-sph-ball}).
Moreover, from (\ref{mv-sph-gen}) we have
$$
\frac{d}{dr}J_v(r)=-\frac{n}{r^{n+1}}\int_{E_r}\left(\frac{\p}{\p
t}-\D\right)v  \, d\mu\, dt
$$
where
$$
J_v(r):=\int_{\p E_r} v\frac{|\nabla H|^2}{\sqrt{|\nabla
H|^2+|H_t|^2}}\, d\tilde A+\frac{1}{r^n}\int_{E_r}R v\, d\mu\, dt.
$$
\end{remark}

When $\Upsilon_{ij}=R_{ij}$, namely we are in the situation of the
Ricci flow, thanks to \cite{P}, there exists a `sub-heat kernel' to
the {\it conjugate heat equation} constructed using the {\it reduced
distance} function discovered by Perelman \cite{P}. For fixed point
$(x_0, 0)$, $\tau=-t$ as above, recall that
$$
L(y, \bar{\tau})=\inf_{\gamma}\mathcal{L}(\gamma)
$$
where the infimum is taken for all arcs $\gamma(\tau)$  joining
$x_0=\gamma(0)$ to $y=\gamma(\bar{\tau})$, and
\begin{equation}\label{L-len}
\mathcal{L}(\gamma)=\int_0^{\bar{\tau}}\sqrt{\tau}(|\gamma'|^2+R(\gamma(\tau),
\tau))\, d\tau
\end{equation}
called the $\mathcal{L}$-length of $\gamma$.
 The `reduced distance' is defined by $\ell(y,
\tau)=\frac{L(y, \tau)}{2\sqrt{\tau}}$. The `sub-heat kernel' is  $
\hat{K}(y, \tau; x_0, 0)=\frac{e^{-\ell(y, \tau)}}{(4\pi
\tau)^{\frac{n}{2}}}. $ One can similarly define the `pseudo heat
ball' $\hat{E}_r$. An important result of \cite{P} asserts that
$$
\left(\frac{\partial}{\partial \tau}-\Delta +R\right)\hat K \le 0.
$$
Restricted in a sufficiently small parabolic neighborhood of $(x_0,
0)$, we know that $\hat{K}(y, \tau; x_0, 0)$ is smooth, and for
small $r$ one can check the compactness of the pseudo-heat ball (cf.
\cite{EKNT}, also the next section). A corollary of Theorem
\ref{mv-sph} is the following result.

\begin{corollary}\label{mv-sph-2}
Let $v\ge 0$ be a smooth function on $\widetilde{M}$.  Let
$\hat{\phi}_r=\hat{K}-r^{-n}.$ Then
\begin{eqnarray}\label{mv-sph2}
v(x_0, 0)&\ge&\int_{\p \hat{E}_r} v\frac{|\nabla
\hat{K}|^2}{\sqrt{|\nabla \hat{K}|^2+|\hat{K}_t|^2}}\, d\tilde
A+\frac{1}{r^n}\int_{\hat{E}_r}R v\, d\mu\,
dt\nonumber\\
&\quad &+\int_{\hat{E}_r}\hat{\phi}_r\left(\frac{\p}{\p
t}-\D\right)v \, d\mu\, dt.
\end{eqnarray}
\end{corollary}
A related new monotonicity quantity is
$$
\hat{J}_v(r)=\int_{\p \hat{E}_r} v\frac{|\nabla
\hat{K}|^2}{\sqrt{|\nabla \hat{K}|^2+|\hat{K}_t|^2}}\, d\tilde
A+\frac{1}{r^n}\int_{\hat{E}_r}R v\, d\mu\, dt.
$$
The virtue of the proof of Corollary \ref{mono2} as well as that of
Theorem \ref{mv-sph} proves the following
\begin{corollary}
\begin{equation}\label{mono-rf1}
\frac{d}{d r} \hat{J}_v(r)\le
-\frac{1}{r^{n+1}}\int_{E_r}\left(\frac{\p}{\p t}-\D\right)v \,
d\mu\, dt.
\end{equation}
\end{corollary}
This, in the case $v=1$ (or more generally any $v$ with
$\left(\frac{\p}{\p t}-\D\right)v\ge0$),   gives a new monotonicity
quanity/formula for the Ricci flow. In \cite{EKNT}, the monotonicity
was proved for the quantity
$$
\hat I_v (r):= \frac{1}{r^n}\int_{\hat{E} _r}\left(|\nabla \log
\hat{K}|^2+R\hat{\psi}_r\right)v\, d\mu dt
$$
where $\hat{\psi}=\log(\hat K r^n)$.  $\hat{I}_v(r)$ is related to
$\hat J_v(r)$ by the relation $ r^n \hat
I_v(r)=n\int_0^r\eta^{n-1}\hat J_v(\eta)\, d\eta. $  The previous
known property that $\hat I_v(r)$ is monotone non-increasing  does
follows from that of $\hat J_v(r)$. In fact, we can rewrite
$$
\hat I_v (r)=\frac{\int_0^r \eta^{n-1} \hat J_v(\eta)\,
d\eta}{\int_0^r \eta^{n-1}\, d\eta}.
$$
Then the result follows from the elementary fact that
$\frac{\int_a^r f(\eta)d\eta}{\int_a^r g(\eta)d\eta}$ is monotone
non-increasing for any $a\in [0, r]$, provided that
$\frac{f(r)}{g(r)}$ is non-increasing. In fact, we have the
monotonicity (non-increasing in both $r$ and $a$) of a more general
quantity
$$
\hat{I}_v(a, r):=\frac{1}{r^n-a^n}\int_{\hat{E}_r\setminus
\hat{E}_a}\left(|\nabla \log \hat{K}|^2+R\hat{\psi}_r\right)v\, d\mu
dt.
$$

 Concerning the spherical mean value theorem, it seems more
natural to consider the fundamental solution of the {\it backward
heat equation} $\frac{\p}{\p t}+\D$ and a solution to the {\it
forward conjugate heat equation} $\frac{\p}{\p t}-\D -R$, since if
$H(y, \tau; x_0, 0)$ is the fundamental solution to $\frac{\p}{\p
\tau}-\D$ and we define the `heat ball', `heat sphere', $\phi_r$ and
$\psi_r$  in the same way as before we can have the following
cleaner result.

\begin{theorem}\label{mv-sph-3}
Let $v$ be a smooth function on $\widetilde{M}$.  Then
\begin{equation}\label{mv-sph3}
v(x_0, 0)=\int_{\p E_r} v\frac{|\nabla H|^2}{\sqrt{|\nabla
H|^2+|H_t|^2}}\, d\tilde A +\int_{E_r}\phi_r\left(\frac{\p}{\p
t}-\D-R\right)v  \, d\mu\, dt.
\end{equation}
\end{theorem}
This looks nicer if $v$ is a solution to the {\it forward conjugate
heat equation}. We also have the  following related result.

\begin{corollary}\label{mono-heat}
Let
$$
J^f_v(r)=\int_{\p E_r} v\frac{|\nabla H|^2}{\sqrt{|\nabla
H|^2+|H_t|^2}}\, d\tilde A, \quad
I^f_v(r)=\frac{1}{r^n}\int_{E_r}|\nabla \log H|^2 v\, d\mu\, dt.
$$
Then
\begin{equation}\label{mono-heat-1}
\frac{d}{d r} J^f_v(r)=-\frac{n}{r^{n+1}}\int_{E_r}\left(\frac{\p
v}{\p t}-\D v-R v\right) \, d\mu\, dt
\end{equation}
and
\begin{equation}\label{mono-heat-2}
\frac{d}{d r}I^f_v(r)=-\frac{n}{r^{n+1}}\int_{E_r} \psi_r
\left(\frac{\p v}{\p t}-\D v-R v\right)  \, d\mu\, dt.
\end{equation}
In particular, $I^f_v(r)$ and $J^f_v(r)$ are monotone non-increasing
if $v$ is a sup-solution to the forward conjugate heat equation.
\end{corollary}

Since most of the above discussion works for any family of metrics
satisfying (\ref{evolution}) we can also apply it to the mean
curvature flow setting. The mean value theorem and related
monotonicity formulae, with respect to the `heat balls',  have been
studied in \cite{E}. Here we just outline the result with respect to
the `heat spheres'. First recall that a family $(M_t)_{t\in(\alpha,
\beta)}$ of a n-dimensional submanifolds of $\mathbb{R}^{n+k}$ moves
by mean curvature if there exist immersions $y_t=y(\cdot, t): M^n\to
\mathbb{R}^{n+k}$ of an $n$-dimensional manifold $M^n$ with images
$M_t=y_t(M^n)$ satisfying the equation
\begin{equation}\label{mcf}
\frac{\partial y}{\partial t}=\vec{H}.
\end{equation}
Here $\vec{H}(p, t)$ denotes the mean curvature vector of $M_t$ at
$y(p, t)$. (We shall still use $H$ to denote the fundamental
solutions.) It was known (see for example \cite{Hu, E}) that the
induced metric $g_{ij}(p, t)=\langle \nabla_i y, \nabla_j y\rangle $
on $M$ satisfying the equation
$$
\frac{\partial}{\partial t}g_{ij}=-2\vec{H}\vec{H}_{ij}
$$
where $\vec{H}_{ij}(p, t)$ is the second fundamental form of $M_t$
at $y(p, t)$. (In this case $R=|\vec{H}|^2$.) By the virtue of
 \cite{H} (see also \cite{E})  we know that the `sup-heat kernel'
\begin{equation}\label{sup-heat}
\bar{K}(y, \tau; x_0, 0)=\frac{1}{(4\pi
\tau)^{\frac{n}{2}}}\exp\left(-\frac{|x_0- y|^2}{4\tau}\right)
\end{equation} viewed as a function on $M$ via the immersion $y_t$,
satisfies that
\begin{equation}
\left(\frac{\partial}{\partial \tau}-\Delta +|\vec{H}|^2\right) \bar
K =\bar{K}\left|\vec{H}-\frac{\nabla^{\perp} \bar K}{\bar
K}\right|^2\ge 0.
\end{equation}
Here  $|x_0- y|$ denotes the Euclidean distance of
$\mathbb{R}^{n+k}$. The proof of Theorem \ref{mv-sph} gives the
following corollary.

\begin{corollary}\label{mv-sph-mcf}
Let $v\ge 0$ be a smooth function on $\mathbb{R}^{n+k}\times(\alpha,
\beta)$. Assume that $M_0$ meets $y_0(x_0)$. Let
$\bar{\phi}_r=\bar{K}-r^{-n}.$ Then
\begin{eqnarray}\label{mv-sph-mcf}
v(x_0, 0)&=&\int_{\p \bar{E}_r} v\frac{|\nabla
\bar{K}|^2}{\sqrt{|\nabla \bar{K}|^2+|\bar{K}_t|^2}}\, d\tilde
A+\frac{1}{r^n}\int_{\bar{E}_r}|\vec{H}|^2 v\, d\mu\,
dt\nonumber\\
&\quad &+\int_{\bar{E}_r}\bar{\phi}_r\left(\frac{\p}{\p
t}-\D\right)v \, d\mu\, dt-\int_{\bar{E}_r}v\bar
K\left|\vec{H}-\frac{\nabla^{\perp} \bar K}{\bar K}\right|^2\, d\mu
dt.
\end{eqnarray}
Here we denote by the same symbols on the functions on
$\mathbb{R}^{n+k}\times(\alpha, \beta)$ and their pull-backs on
$M\times(\alpha, \beta)$.
\end{corollary}
Integrating (\ref{mv-sph-mcf}) one can recover the mean value
formula in \cite{E},  noticing, along with (\ref{equ1}),   that for
any function $f$ (regular enough)
$$
\frac{n}{r^n}\int_0^r \eta^{n-1}\int_{\bar{E}_\eta} \bar K f\, d\mu
\, d t\, d\eta =\int_0^r \frac{n}{\eta^{n+1}}\int_{\bar E_\eta} f\,
d\mu\, dt\, d\eta.
$$

If we let
$$
\bar J_v(r)=\int_{\p \bar{E}_r} v\frac{|\nabla
\bar{K}|^2}{\sqrt{|\nabla \bar{K}|^2+|\bar{K}_t|^2}}\, d\tilde
A+\frac{1}{r^n}\int_{\bar{E}_r}|\vec{H}|^2 v\, d\mu\, dt
$$
we have that
\begin{equation}\label{mono-sph-mcf}
\frac{d}{d r}\bar J_v(r)=\frac{n}{r^{n+1}}\left(-\int_{\bar E
_r}\left(\frac{\p}{\p t}-\D\right)v\, d\mu\, dt +\int_{\partial \bar
E _r}  \frac{v\bar K}{|\tilde \nabla \bar
K|}\left|\vec{H}-\frac{\nabla^{\perp} \bar K}{\bar K}\right|^2\,
d\tilde A \right).
\end{equation}
Here $\tilde{\nabla}\bar{K}=\langle \nabla \bar{K}, \frac{\partial
\bar{K}}{\partial t}\rangle$. The equation (\ref{mono-sph-mcf})
gives a new monotonicity formula in case that $v$ is a nonnegative
sub-solution to the heat equation. Note that it was previously
proved in \cite{E} that
$$
\bar{I}(r)=\frac{1}{r^n}\int_{\bar E_{r}}\left(|\nabla \log \bar
K|^2+|\vec{H}|^2 \psi_r\right)\, d\mu\, dt
$$
is monotone non-decreasing. As the above Ricci flow case, the
monotonicity on $\bar{I}(r)$ follows from the monotonicity of $\bar
J_v(r)$ for $v\equiv 1$. In fact we also have that
$$
\bar{I}(a, r)=\frac{1}{r^n-a^n}\int_{\bar E_{r}\setminus
\bar{E}_{a}}\left(|\nabla \log \bar K|^2+|\vec{H}|^2 \psi_r\right)\,
d\mu\, dt
$$
is monotone non-decreasing in both $r$ and $a$.

\section{Local regularity theorems}
The monotonicity formula proved in \cite{EKNT} (as well as the one
proved in Proposition 5.4 of \cite{N3}), together with Hamilton's
compactness theorem \cite{H2} and  the arguments in Section 10 of
\cite{P}, allows us to formulate a $\epsilon$-regularity theorem for
Ricci flow Theorem \ref{Theorem 3.1} (Theorem \ref{Theorem 3.2}).
The closest results of this sort are the ones for mean curvature
flow \cite{E2, Wh}. The result here is influenced by \cite{E2}.

We first recall some elementary properties of the so-called `reduced
geometry' (namely the geometry related to the functional
$\mathcal{L}(\gamma)$). Let $(M, g(t))$ be a solution to Ricci flow
on $M\times [0, T]$. Let $(x_0, t_0)$ be a fixed point with
$\frac{T}{2}\le t_0\le T$. We can define the $\mathcal{L}$-length
functional for any path originated from $(x_0, t_0)$. Let
$\ell^{(x_0, t_0)}(y, \tau)$ be the `reduced distance' with respect
to $(x_0, t_0)$.

Here we mainly follow \cite{N3} (Appendix), where the simpler case
for a fixed Riemannian metric was detailed. It seems more natural to
consider $\tM =M \times [0, t_0]$. Perelman defined the ${\mathcal
L}$-exponential map  as follows. First one
 defines the ${\mathcal L}$-geodesic to be the critical point of the
above ${\mathcal L}$-length functional in (\ref{L-len}). The
geodesic equation
$$
\nabla_X X -\frac{1}{2}\nabla R +\frac{1}{2\tau}X+2Ric(X, \cdot)=0,
$$
where $X=\frac{d \gamma}{d\tau}$,  is derived in \cite{P},  which
has (regular) singularity at $\tau=0$. It is desirable to change the
variable to $\sigma=2\sqrt{\tau}$. Now we can define
$\bar{g}(\sigma)=g(\tau)=g(\frac{\sigma^2}{4})$. It is easy to see
that $\frac{\partial}{\partial \sigma}\bar{g}=\sigma Ric$. The
$\mathcal{L}$-length has the form
$$
\mathcal{L}(\gamma)=\int_0^{\bar{\sigma}}\left(\left|\frac{d
\gamma}{d \sigma}\right|^2+\frac{\sigma^2}{4}R\right)\, d\sigma.
$$
The $\mathcal{L}$-geodesic has the form
$$
\nabla_{\bar{X}}\bar{X}+\sigma Ric(\bar X,
\cdot)-\frac{\sigma^2}{8}\nabla R=0
$$
where $\bar{X}=\frac{d \gamma}{d\sigma}$. By the theory of ODE we
know that for any $v\in T_{x_0}M$ there exists a ${\mathcal
L}$-geodesic $\gamma(\sigma)$ such $\frac{ d}{d
\sigma}\left(\gamma(\sigma)\right)|_{\sigma =0}=v$, where
$\sigma=2\sqrt{\tau}$. Then as in \cite{P} one defines the
${\mathcal L}$-exponential map as
$$
{\mathcal L}\exp_{v}(\sigmab):=\gamma_{v}(\sigmab)
$$
where $\gamma_v(\sigma)$ is a ${\mathcal L}$-geodesic satisfying
that
$
\lim_{\sigma \to 0}\frac{d
}{d\sigma}\left(\gamma_v(\sigma)\right)=v.
$
The space-time exponential map $\texp: T_{x_0}M\times [0,
2\sqrt{t_0}] \to \tM $ is defined as
  $$\texp(\tilde v^a)=({\mathcal
L}\exp_{\frac{v^a}{a}}(a), a),$$
 where $\tilde v^a=(v^a, a)$. (Here we abuse the notation $\tM$ since in terms of
$\sigma$, $\tM=M\times [0, 2\sqrt{t_0}]$.) Denote $\frac{v^a}{a}$
simply by $v^1$ and $(v^1, 1)$ by $\tilde v^1$. Also let $\tilde
\gamma_{\tilde v^a}(\eta)=\texp(\eta\tilde v^a)$. It is easy to see
that
\begin{equation}\label{geo}
\tilde \gamma_{\tilde v^a}(\eta)=\tilde \gamma_{\tilde v^1}(\eta a).
\end{equation}
This implies the following lemma.
\begin{lemma}\label{exp}
$$
d\texp|_{(0, 0)}=\text{identity}.
$$
In particular, there exists a neighborhood $U\subset
T_{x_0}M\times[0, 2\sqrt{t_0}]$ near $(x_0, 0)$ so that restricted
on $U\cap \left( T_{x_0}M\times(0, 2\sqrt{t_0}]\right)$, $\texp$ is
a diffeomorphism.
\end{lemma}
\begin{proof} The first part follows from taking derivative on
(\ref{geo}).  Tracing the proof of the implicit function theorem
(see for example \cite{Ho}) gives the second part.
\end{proof}

  Based on this fact we can show the following result analogous to the Guass
lemma.
  \begin{corollary}\label{guass}The
${\mathcal L}$-geodesic $\gamma_{v}(\sigma)$ from $(x_0, 0)$ to $(y,
\sigmab)$ is ${\mathcal L}$-length minimizing for $(y, \sigmab)$
close to $(x_0, 0)$, with respect to metric $\tilde g(x, \sigma)$ on
$\tM$, where $\tilde{g}(x, \sigma)=\bar{g}(x, \sigma)+d\sigma^2$ (or
any other compatible topology).
\end{corollary}
\begin{proof} We adapt the parameter $\tau$ for this matter. By the
above lemma we can define function $\underline{\mathcal{L}}(y,
\tau)$ to be the length of the $\mathcal{L}$-geodesic jointing from
$(x_0, 0)$ to $(y, \tau)$. For $(y, \tau)$ close to $(x_0, 0)$ it is
a smooth function. It can be shown by simple calculation as  in
\cite{P} that
\begin{equation}\label{per1}
\nabla \underline{\mathcal{L}}(y, \tau)=2\sqrt{\tau_1}X
\end{equation}
where $\nabla$ is the spatial gradient with respect to $g(\tau)$,
$X=\gamma'(\tau)$ with $\gamma(\tau)$ being the
$\mathcal{L}$-geodesic joining $(x_0, 0)$ to $(y, \tau)$, and
\begin{equation}\label{per2}
\frac{\partial}{\partial \tau}\underline{\mathcal{L}}(y,
\tau)=\sqrt{\tau}(R-|X|^2).
\end{equation}
Now for any curve $(\eta(\tau), \tau)$ joining $(x_0, 0)$ to some
$(y, \tau_1)$, we have that
\begin{eqnarray*}
\underline{\mathcal{L}}(\eta(\tau_1),
\tau_1)&=&\int_0^{\tau_1}\frac{d\underline{\mathcal{L}}}{d\tau}\, d\tau \\
&=& \int_0^{\tau_1}\langle \nabla \underline{\mathcal{L}},
\eta'\rangle
+\frac{\partial \underline{\mathcal{L}}}{\partial \tau}\, d\tau\\
&=& \int_0^{\tau_1}2\sqrt{\tau}\langle X, \eta'\rangle
+\sqrt{\tau}(R-|X|^2)\, d\tau\\
&\le & \int_0^{\tau_1}\tau\left(R+|\eta'|^2\right)\, d\tau\\
&=&\mathcal{L}(\eta).
\end{eqnarray*}
This shows that the $\mathcal{L}$-geodesic is the
$\mathcal{L}$-length minimizing path joining $(x_0, 0)$ with
$(y,\tau_1)$.
\end{proof}

Lemma \ref{exp} particularly implies that ${\mathcal
L}\exp_{v}(\tau)$ (which is defined as $\gamma_v(\tau)$ with
$v=\lim_{\tau\to 0}2\sqrt{\tau}\gamma'(\tau)$) is a diffeomorphism
if $|v|\le 1$, for $\tau\le \epsilon$, provided that $\epsilon$ is
small enough. Here we identify $T_{x_0}M$ with the Euclidean space
using the metric $g(x_0, t_0)$.
 One can similarly define the cut point (locus) and the first
conjugate point (locus) with respect to the ${\mathcal L}$-length
functional. For example, $y=\gamma_v(\sigmab)$ (with
$\sigmab<2\sqrt{t_0}$) is a cut point if the ${\mathcal L}$-geodesic
$\gamma_{v}(\sigma)$  is minimizing up to $\sigmab$ and fails to be
so for any $\sigma>\sigmab$ (the case that $\sigmab=2\sqrt{t_0}$
needs a different definition which shall be addressed later).
Similarly one defines the first ${\mathcal L}$-conjugate point. Now
we  let
$$
D(\sigmab)\subseteq T_{x_0}M
$$to be  the collection of vectors $v$ such that $({\mathcal L}\exp_v(\sigma),
\sigma)$ is a ${\mathcal L}$-geodesic along which there is no
conjugate for all  $\sigma<\sigmab$. Similarly we let
$$
\Sigma(\sigmab)\subseteq T_{x_0}M
$$
to be the collection of vectors $v$ such that $({\mathcal
L}_v(\sigma), \sigma)$ is a minimizing ${\mathcal L}$-geodesic up to
$\sigmab$. One can see easily that $\Sigma(\sigma)\subset
D(\sigma)$, and both
 $D(\sigma)$ and $\Sigma(\sigma)$ decrease (as sets) as $\sigma$
increases. For any measurable subset $A\subset T_{x_0}M$ we can
define
$$
D_A(\sigma)=A\cap D(\sigma) \quad \text{and}\quad
\Sigma_A(\sigma)=A\cap \Sigma(\sigma).
$$
 It
is exactly the same argument as the classical case to show that for
a cut point $(y, \bar \sigma)$ with $\sigmab<2\sqrt{t_0}$, either it
is a conjugate point (namely $(y, \sigmab)$ is a critical value of
$\texp$) or there are two minimizing ${\mathcal L}$-geodesics
$\gamma_{v_1}(\sigma)$ and $\gamma_{v_2}(\sigma)$ joining $x_0$ with
$y$ at $\sigmab$. For $\sigmab=2\sqrt{t_0}$, we can use this
property to define the cut points. Namely $(y, 2\sqrt{t_0})$ is
called a cut point if  either it is a critical value of
$\widetilde{\exp}$ or there exists two minimizing
$\mathcal{L}$-geodesics joining $(x_0, 0)$ with $(y, 2\sqrt{t_0})$.
For a fixed $v\in T_{x_0}M$ we can define $\bar{\sigma}(v)$ to be
the first $\bar{\sigma}$, if it is smaller than $2\sqrt{t_0}$, such
that $\gamma_v(\sigma)$ is $\mathcal{L}$-minimizing for all
$\sigma\le\bar{\sigma}$ and no longer so or not defined  for
$\sigma>\bar{\sigma}$. If the original solution to the Ricci flow
$g(t)$ is an ancient solution, then $\widetilde{M}=M\times[0,
\infty)$. We define $\bar{\sigma}=\infty$ if $\gamma_v(\sigma)$ is
always $\mathcal{L}$-minimizing, in which case we call
$\widetilde{\gamma}_{\tilde v^{1}}(\sigma)$, with $\tilde{v}^1=(v,
1)$, a $\mathcal{L}$-geodesic ray in $\widetilde{M}$. For the case
$\widetilde{M}=M \times [0, 2\sqrt{t_0}]$ with finite $t_0$, define
$\bar{\sigma}(v)=2\sqrt{t_0}+$ (slightly bigger than $2\sqrt{t_0}$)
if $\gamma_v(\sigma)$ is minimizing up to $2\sqrt{t_0}$ and it is
the unique $\mathcal{L}$-geodesic joining to
$(\gamma_v(2\sqrt{t_0}), 2\sqrt{t_0})$ from $(x_0, 0)$. One can show
that $\bar{\sigma}(v)$ is a continuous function on $T_{x_0}M$.
Define
$$
\widetilde{\Sigma}=\{(\sigma v, \sigma)\, |\, v \in T_{x_0}M,  0\le
\sigma< \bar{\sigma}(v)\}.
$$
It can also be shown that $\widetilde{\Sigma}$ is open (relatively
in $T_{x_0}M\times [0, 2\sqrt{t_0}]$) and can be checked that
$\widetilde{\exp}|_{\widetilde{\Sigma}}$ is a diffeomorphism. The
cut locus (in the tangent half space $T_{x_0}M \times \mathbb{R}_+$,
or more precisely $T_{x_0}M\times [0, 2\sqrt{t_0}]$)
$$\mathcal{C}= \{ (\bar{\sigma}(v)v, \bar{\sigma}(v))\in T_{x_0}M
\times[0, 2\sqrt{t_0}]\subset T_{x_0}M\times \mathbb{R}_+ \}.$$ It
is easy to show, by  Fubini's theorem, that $\mathcal{C}$ has zero
($n+1$-dimensional) measure. As in the classical case
$\widetilde{\exp}(\widetilde{\Sigma} \cup
\mathcal{C})=\widetilde{M}$. Here comparing with the classical
Riemannian geometry, $T_{x_0}M$ which can be identified with
$T_{x_0}M\times\{1\}\subset T_{x_0}M \times \mathbb{R}_+$, plays the
same role as the unit sphere of the tangent space in the Riemannian
geometry.

However, there are finer properties on the cut locus. Let
$C(\sigma_0)$ be the cut points at $\sigma_0$ and ${\mathcal
C}(\sigma_0)$ the corresponding vectors in $T_{x_0}M$. Namely
$\mathcal{C}(\sigma_0)=\{v \, | \bar{\sigma}(v)=\sigma_0\}$ and
$C(\sigma_0)=\mathcal{L}\exp_{\mathcal{C}(\sigma_0)}(\sigma_0).$
Then
$\widetilde{\exp}\left(\left(\widetilde{\Sigma}\cap(T_{x_0}M\times\{\sigma\})\right)\cup
\left(\sigma\mathcal{C}(\sigma)\times\{\sigma\}\right)\right)=M
\times\{\sigma\}$. More importantly, for any $\sigma<2\sqrt{t_0}$ ,
$C(\sigma)\times\{ \sigma\}\subseteq M\times \{\sigma\}$ has
($n$-dimensional) measure zero, provided that the metrics $g(x,
\sigma)$ are sufficiently smooth  on $\widetilde{M}$.  By Sard's
theorem, we know that the set of the conjugate points in $C(\sigma)$
has zero measure. Hence to prove our claim we only need to show the
second type cut points (the ones having more than one minimizing
$\mathcal{L}$-geodesics ending at) in $C(\sigma)$ has measure zero
in $M\times \{\sigma\}$. These exactly are the points on which
$\nabla L(\cdot, \sigma)$ is not well-defined. (Namely $L(\cdot,
\sigma)$ fails to be differentiable.) Let $\bar{L}=\sigma L$. It was
proved in \cite{Ye} (see also \cite{CLN}) that $\bar{L}(y, \sigma)$
is locally Lipschitz on $M\times [0, \sigmab]$. In particular, for
every $\sigma>0$, $\ell(y, \sigma)$ is locally Lipschitz as a
function of $y$.  The claimed result that $C(\sigma)$ has measure
zero follows from Rademacher's theorem on the almost everywhere
smoothness of a Lipschitz function. Once we have this fact, if
defining that
$$\widetilde{\Sigma}(\sigma)=\{v\, |\, (\sigma v, \sigma)\in
\widetilde{\Sigma}\},$$ for any integrable function $f(y, \sigma)$
we have that
$$
\int_M f(y, \sigma)\,
d\mu_{\sigma}=\int_{\widetilde{\Sigma}(\sigma)} f(y, \sigma)\, J(v,
\sigma)\, d\mu_{euc}(v)
$$
where $y$ is a function of $v$ through the relation
$v=(\mathcal{L}\exp_{(\cdot)} (\sigma))^{-1}(y)$, $J(v, \sigma)$ is
the Jacobian of $\mathcal{L}\exp_{(\cdot)}(\sigma)$ and $d\mu_{euc}$
is the volume form of $T_{x_0}M=\mathbb{R}^n$ (via the metric of
$g(t_0)$). It is clear that $\widetilde{\Sigma}(\sigma)\subset
\Sigma(\sigma)$ and $\widetilde{\Sigma}(\sigma)$ decreases as
$\sigma$ increases.

 The above discussion can be translated in terms of $\tau$.
 Let
 \begin{equation}\label{sub-heat}\hat K(y, \tau; x_0, t_0)=\frac{1}{(4\pi \tau)^{\frac{n}{2}}}\exp\left(-\ell^{(x_0,
t_0)}(y, \tau)\right).\end{equation}It was proved in \cite{P} that
if $y$ lies in $\mathcal{L}\exp_{\widetilde{\Sigma}(\tau)}(\tau)$,
$$
\frac{d}{d\tau} \left(\hat{K}(y, \tau; x_0, t_0)J(v, \tau)\right)
\le 0
$$
where $y=\mathcal{L}\exp_{v}(\tau)$, which further implies that the
`reduced volume'
$$
\theta(\tau):=\int_M \hat K(y, \tau; x_0, t_0)\, d\mu_\tau
$$
is monotone non-increasing in $\tau$.
 Let $\hat{E}_r$ be the `pseudo
 heat ball' and
$$
\hat{I}^{(x_0, t_0)}(r):= \frac{1}{r^n}\int_{\hat{E} _r}\left(
|\nabla \log \hat{K}|^2+R\hat{\psi}_r\right) \, d\mu dt
$$
where $\hat{\psi}_r=\log(\hat{K} r^n)$.  The following result, which
can be viewed as a localized version of Perelman's above
monotonicity of $\theta(\tau)$, was proved in \cite{EKNT}. (This
also follows from Corollary \ref{mv-sph-2}).

\begin{theorem}\label{eknt} a) Let $(M, g(t))$ be a solution to the Ricci flow
as above. Let $\ell^{(x_0, t_0)}(y,\tau)$, $\hat{K}$ and
$\hat{I}^{(x_0, t_0)}(r)$ be defined as above. Then $\frac{d}{
dr}\hat{I}^{(x_0, t_0)}(r)\le 0$. If the equality holds for some
$r$, it implies that $(M, g(t))$ is a gradient shrinking Ricci
soliton on $\hat{E}_r$.

b) If $(M, g)$ be a complete Riemannian manifold with $Ric\ge
-(n-1)k^2 g$. Let $\overline {H}(x, y, \tau)=\overline{H}(\bar{d}(x,
y), \tau)$ be the fundamental solution  of backward heat equation
$\frac{\partial}{\partial \tau}-\Delta$ on the space form $\bar{M}$
with constant curvature $-k^2$, where $\bar{d}(\cdot, \cdot)$ is the
distance function of $\bar{M}$. Define $\hat{I}^{(x_0, t_0)}(r)$
similarly using the `sub heat kernel' $\hat{H}(y, \tau; x_0,
0)=\overline{H}(d(x, y), \tau)$. Then $\frac{d}{d r}\hat{I}^{(x_0,
t_0)}(r) \le 0$. The equality for some $r$  implies that $(M, g)$ is
isometric to $\bar{M}$ on $\hat{E}_r$.
\end{theorem}

In \cite{EKNT}, it was shown that  for a solution $g(x, t)$ to Ricci
flow defined on $M\times [0, T]$, $\hat{I}^{(x, t)}(r, g)$ is
well-defined if $t\ge \frac{3}{4}T$ and $r$ is sufficiently small.
We adapt the notation of \cite{P} by denoting the parabolic
neighborhood $B_t(x, r)\times [t+\Delta t, t]$ (in the case $\Delta
t<0$) of $(x, t)$ by $P(x, t, r, \Delta t)$. Let $\widetilde{U}=
U\times[t_1, t_0]$ be a parabolic neighborhood of $(x_0, t_0)$. Let
$r_0$ and $\rho_0$ are sufficient small so that $P(x_0, t_0, r_0,
-r_0^2)\subset \widetilde{U}$ and $\hat{I}^{(x, t)}(\rho)$ is
defined for all $\rho\le \rho_0$ and $(x, t)\in P(x_0, t_0, r_0,
-r_0^2)$. Motivated by Perelman's pseudo-locality theorem, we prove
the following local curvature estimate result for the Ricci flow.

\begin{theorem}\label{Theorem 3.1} There exist positive  $\epsilon_0$ and $C_0$ such that that if
$g(t)$ is a solution of Ricci flow on a parabolic neighborhood
$\widetilde{U}$ of $(x_0, t_0)$ and  if
\begin{equation}\label{3.1}\hat{I}^{(x, t)}(\rho_0)\ge 1-\e_0,
\end{equation}
for $(x,t)\in P(x_0, t_0, r_0, -r_0^2)$
 then
\begin{equation}\label{3.2}
|Rm|(y, t)\le \frac{C_0}{r_0^2} \end{equation} for any $(y, t)\in
P(x_0, t_0, \frac{1}{2}r_0, -\frac{1}{4}r_0^2)$.
\end{theorem}
\begin{proof} We prove the result by contradiction argument via Hamilton's compactness theorem \cite{H2}.
Notice that both the assumptions and the conclusion are scaling
invariant. So we can assume $r_0=1$ without the loss of the
generality. Assume the conclusion is not true. We then have
 a sequence of counter-examples to the theorem. Namely there exist
 $(\widetilde{U}_j, g_j(t))$, parabolic neighborhood of   $(x^j_0, t^j_0)$
and  $P^j_1=P(x^j_0, t^j_0, 1, -1)$ with
$$
\hat{I}^{(x, t)}(\rho_0, g_j)\ge 1-\frac{1}{j}
$$
for all $(x, t)\in P^j_1$ (here we write $\hat I$ explicitly on its
dependence of the metric) but
$$
Q_j:= \sup_{P^j_{\frac{1}{2}}} |Rm(g_j)|(y, t)\ge j\to \infty
$$
where $P^j_{\frac{1}{2}}=P^j(x^j_0, t^j_0, \frac{1}{2},
-\frac{1}{4})$.

By the argument in Section 10 of \cite{P}  we can  find $(\bar{x}_j,
\bar{t}_j)\in P^j_{\frac{3}{4}}$ such that
 $|Rm|(y, t)\le 4\bar Q_j$,
for any $(y, t)\in  \bar{P}_j:=P(\bar{x}_j, \bar{t}_j, H_j
\bar{Q}_j^{-\frac{1}{2}}, -H^2_j\bar{Q}_j^{-1})$, with $H_j\to
\infty$ as $j\to \infty$. Here $\bar{Q}_j=|Rm|(\bar{x}_j,
\bar{t}_j)$. This process can be done through two steps as in Claim
1 and 2 of  Section 10 of \cite{P}. Let $\tilde g_j(t):=\bar{Q}_j
g(\bar{t}_j +\frac{t}{\bar{Q}_j})$. Now we consider two cases.

{\it Case 1}: The injectivity radius of $(\widetilde{U}_j, \tilde
g_j(0))$ is bounded from below uniformly at $\bar{x}_j$. In this
case by Hamilton's compactness theorem \cite{H2} we can conclude
that $(\bar{P}_j, \bar{x}_j, \tilde{g}_j)$ converges to $(M_\infty,
x_\infty, g_\infty)$, an ancient solution to Ricci flow with
$|Rm|(x_\infty, 0)=1$. On the other hand, for any $\rho\le \rho_0$,
by the monotonicity
$$
1\ge \hat{I}^{(\bar{x}_j, \bar{t}_j)}(\rho, g_j)\ge 1-\frac{1}{j}
$$
which then implies that for any $\rho\le \rho_0 \bar{Q}_j$,
$$
1\ge \hat{I}^{(\bar{x}_j, 0)}(\rho, \tilde g_j)\ge 1-\frac{1}{j}.
$$
This would imply that for any $\rho>0$,
$$
\hat{I}^{(x_\infty, 0)}(\rho, g_\infty)\equiv 1.
$$
By Theorem \ref{eknt}, it then implies that there exists a function
$f$ such that
$$
R_{ij}+f_{ij}+\frac{1}{2t}g_{ij}=0.
$$
Namely, $(M_\infty, g_\infty)$ is a non-flat gradient shrinking
soliton, which become singular at $t=0$. This is a contradiction to
$|Rm|(x_\infty, 0)=1$.

{\it Case 2}: The injectivity radius of  $(\widetilde{U}_j, \tilde
g_j(0))$ is not bounded from below uniformly at $\bar{x}_j$. By
passing to a subsequence, we can assume that  the injectivity radius
of  $(\widetilde{U}_j, \tilde g_j(0))$  at $\bar{x}_j$, which we
denote by $\lambda_j$, goes to zero. In this case we re-scale the
metric $\tilde g_j(t)$ further by letting $\tilde
g^*_j(t)=\frac{1}{\lambda_j^2}\tilde g_j(\lambda_j ^2 t)$. The new
sequence $(\bar{P}_j, \bar{x}_j, \tilde g^*_j)$  will have the
required injectivity radius lower bound, therefore converges to a
flat limit $(M^*_\infty, x^*_\infty, g^*_\infty(t))$ with the
injectivity radius at $x^*_\infty$ being equal to one. On the other
hand, by the monotonicity of $\hat{I}^{(\bar{x}_j, 0)}(\rho, \tilde
g^*_j)$ and the similar argument as in {\it Case 1} we can conclude
 that $\hat{I}^{(x^*_\infty, 0)}(\rho, g^*_\infty)\equiv 1$ for any $\rho>0$.
 Therefore, by Theorem \ref{eknt},
 we conclude that that
$(M^*_\infty, g^*_\infty)$ must be isometric to $\R^n$. This is a
contradiction!
\end{proof}

If $t_0$ is singular time we may define the `density function' at
$(x_0, t_0)$ by
$$
\hat I^{(x_0, t_0)}(\rho):=\liminf_{(x, t)\to (x_0, t_0)} \hat
I^{(x, t)}(\rho).
$$
Then we may conclude that $(x_0, t_0)$ is a smooth point if $\hat
I^{(x_0, t_0)}(\rho)\ge 1-\epsilon_0$ for some $\rho_0>0$.

 Same result can be formulated for the
localized entropy. Recall that in \cite{P}, for any $(x_0, t_0)$ one
can look at $u$ the fundamental solution to the {\it backward
conjugate heat equation} $\left(\frac{\p}{\p \tau}-\D
+R\right)u(x,\tau)=0$ centered at $(x_0, t_0)$, one have that
$$
\left(\frac{\p}{\p \tau}-\D
+R\right)(-v)=2\tau|R_{ij}+f_{ij}-\frac{1}{2\tau}|^2u
$$
where $\tau=t_0-t$, $v=\left[\tau \left(2\D f -|\nabla
f|^2+R\right)+f-n\right]u$ and $f$ is defined by
$u=\frac{e^{-f}}{\left(4\pi \tau\right)^{\frac{n}{2}}}$. In
particular, when $M$ is compact one has the entropy formula
$$
\frac{d}{ d t}\int_M (-v)\, d\mu_t =-2\int_M\
\tau|R_{ij}+f_{ij}-\frac{1}{2\tau}|^2u\, d\mu_t.
$$

 In \cite{N3} we observed
that the above entropy formula can be localized. In deed, let $\bar
L^{(x_1, t_1)}(y, t)=4(t_1-t)\ell^{(x_1, t_1)}(y, t_1-t)$ and let
$$
\psi^{(x_1, t_1)}_{t_2, \rho}(x, t):=\left(1-\frac{\bar{L}^{(x_1,
t_1)}(x,t)+2n(t-t_2)}{\rho^2}\right)_{+}.
$$
It is easy to see that
$$
\heat \psi^{(x_1, t_1)}_{t_2, \rho}(x, t)\le 0.
$$
The Proposition 5.5 of \cite{N3} asserts that
\begin{equation}\label{3.3}
\frac{d}{d t}\left(\int_M -v\psi^{(x_1, t_1)}_{t_2, \rho}\,
d\mu_t\right)\le -2\int_M \tau
\left(|R_{ij}+\nabla_{i}\nabla_{j}f-\frac{1}{2\tau}g_{ij}|^2\right)
u\psi^{(x_1, t_1)}_{t_2, \rho}\, d\mu_t .
\end{equation}
One can then define
$$
\nu^{(x_0, t_0)}(\rho, \tau,  g)=\int_M -v(y, \tau)\psi^{(x_0,
t_0)}_{t_0, \rho}(y, \tau)\, d\mu_\tau(y).
$$
The similar argument as in the  Section 10 of \cite{P} (or the above
proof of  Theorem \ref{Theorem 3.1}) shows the following result.
First we assume that $g(t)$ is a solution of Ricci flow on $U\times
[0, T)$. Fix a space-time point $(x_0, t_0)$ with $t_0< T$. Let
$\rho_0>0$ $\tau_0$ and $r_0$ be positive constants which are small
enough such that $\nu^{(x, t)}(\rho, \tau, g)$, with $\rho
\le\rho_0$ and $\tau\le \tau_0$,  is defined for all $(x, t)$ in
$P(x_0, t_0, r_0, -r_0^2)$.

\begin{theorem}\label{Theorem 3.2}  There exist positive constants
$\e_0$ and $C_0$ such that for $(U, g(t))$ a solution to the Ricci
flow, $(x_0, t_0)$ and $P(x_0, t_0, r_0, -r_0^2)$ as above, if
\begin{equation}\label{3.4}\nu^{(x, t)}(\rho_0, \tau_0, g)\le \e_0,
\end{equation} for $(x,t)\in P(x_0, t_0, r_0, -r_0^2)$ then
\begin{equation}\label{3.5}
|Rm|(y, t)\le \frac{C_0}{r_0^2} \end{equation} for any $(y, t)\in
P(x_0, t_0, \frac{1}{2}r_0, -\frac{1}{4}r_0^2)$.
\end{theorem}

In the proof one needs to replace Theorem \ref{eknt} by a similar
rigidity type result  (Corollary 1.3 of \cite{N2} on page 371),
formulated in terms of  the entropy,  to obtain the contradiction
for the collapsing case.

In the view of the works \cite{E2,  Wh}, one can think both
$1-\hat{I}^{(x,t)}(\rho, g)$ and $\nu^{(x, t)}(\rho,\tau, g)$ as
local `densities' for Ricci flow.

\section{Further discussions}

First we show briefly that if $\hat{J}(r_1)=\hat{J}(r_2)$ for some
$r_2>r_1$, it implies that $g(t)$ is a gradient shrinking soliton in
$\hat{E}_{r_2}\setminus \hat{E}_{r_1}$. By the proof of Theorem
\ref{mv-sph}, the equality implies that on $\hat{E}_{r_2}\setminus
\hat{E}_{r_1}$, $(\frac{\partial}{\partial \tau} -\Delta +R)\hat
K=0$, which is equivalent to
\begin{equation}\label{heat-ell}
\ell_\tau -\Delta \ell +|\nabla \ell|^2-R +\frac{n}{2\tau}=0.
\end{equation}
By Proposition 9.1 of \cite{P} we have that
\begin{equation}\label{entropy-pw}
\left(\frac{\partial}{\partial \tau} -\Delta
+R\right)v=-2\tau\left|R_{ij}+\ell_{ij}-\frac{1}{2\tau}g_{ij}\right|^2\hat{K}
\end{equation}
where
$$
v=\left(\tau(2\Delta \ell -|\nabla \ell|^2+R)+\ell-n\right)\hat{K}.
$$
On the other hand, (7.5) and (7.6) of \cite{P} implies that $\ell$
satisfies the first order PDE:
\begin{equation}\label{hyper-ell}
-2\ell_\tau -|\nabla \ell|^2+R-\frac{1}{\tau}\ell =0.
\end{equation}
This together with (\ref{heat-ell}) implies that $v=0$ on
$\hat{E}_{r_2}\setminus \hat{E}_{r_1}$. The result now follows from
(\ref{entropy-pw}).

 The space-time divergence theorem from  Section 3,
Lemma \ref{divergence},  allows us to write monotonic quantities
$J_v(r)$, $I_v(r)$ ($\hat{J}_v(r)$, $\hat{I}_v(r)$ as well), defined
in Section 3 and 4, in a nicer form when $v=1$.

Applying the space time divergence theorem to $\tilde
X=\frac{\partial}{\partial t}$ we have that
\begin{equation}\label{div-c1}
\int_{E_r} R\, d\mu\, dt =-\int_{\partial E_r} \langle
\frac{\partial}{\partial t}, \tilde{\nu}\rangle \,d \tilde
A=\int_{\partial E_r} \frac{H_t}{\sqrt{|\nabla H|^2+|H_t|^2}}\,
d\tilde A.
\end{equation}
Therefore we have that ($J(r):=J_{1}(r)$)
\begin{equation}\label{J-new}
J(r)=\int_{\partial E_r}\frac{|\nabla H|^2-HH_{\tau}}{\sqrt{|\nabla
H|^2+|H_{\tau}|^2}}\, d\tilde A.
\end{equation}
Now the quantity is expressed solely in terms of  the surface
integral on $\partial E_r$. For the Ricci flow case we then have
that
\begin{equation}\label{Jhat-new}
\hat{J}(r)=\int_{\partial \hat{E}_r}\frac{|\nabla
\hat{K}|^2-\hat{K}\hat{K}_\tau}{\sqrt{|\nabla
\hat{K}|^2+|\hat{K}_{\tau}|^2}}\, d\tilde A \end{equation} is
monotone non-increasing, where $\hat{K}$ is the `sub-heat kernel'
defined in (\ref{sub-heat}). For the mean curvature flow case we
have that
\begin{equation}\label{Jbar-new}
\bar{J}(r)=\int_{\partial \bar{E}_r}\frac{|\nabla
\bar{K}|^2-\bar{K}\bar{K}_\tau}{\sqrt{|\nabla
\bar{K}|^2+|\bar{K}_\tau|^2}}\, d\tilde A
\end{equation}
is monotone non-decreasing, where $\bar{K}$ is `sup-heat kernel'
defined by (\ref{sup-heat}). Notice that the numerator $|\nabla
H|^2-H H_\tau$, after dividing $H^2$,  is the expression of Li-Yau
in their celebrated gradient estimate \cite{LY}.

Applying Lemma \ref{divergence} to the vector field $\tilde X=\psi_r
\frac{\partial }{\partial t}$ we have that
\begin{equation}\label{div-c2}
\int_{E_r}R \psi_r \, d\mu\, dt =\int_{\partial E_r}\frac{\partial
\psi_r}{\partial t}\, d\mu\, dt =\int_{\partial E_r}\frac{\partial
}{\partial t}\log H\, d\mu\, dt.
\end{equation}
Here we have used that $\psi_r=0$ on $\partial E_r$ and
$$
\lim_{s\to 0}\int_{P^2_s} \psi_r\, d\mu_s =0.
$$
This was also observed in \cite{EKNT}. Hence
\begin{eqnarray}\label{I-new}
I(r)&=&\frac{1}{r^n}\int_{E_r}\left(|\nabla \log H|^2+R
\psi_r\right)\, d\mu\, dt\\
&=&\frac{1}{r^n}\int_{E_r}\left(|\nabla \log H|^2-(\log
H)_\tau\right)\, d\mu\, dt.\nonumber
\end{eqnarray}
In particular, for the Ricci flow we have the non-increasing
monotonicity of
\begin{equation}\label{Ihat-new}
\hat{I}(r)=\frac{1}{r^n}\int_{\hat{E}_r}\left(|\nabla \log
\hat{K}|^2-(\log \hat{K})_\tau\right)\, d\mu\, dt
\end{equation}
and for the mean curvature flow there exists the monotone
non-decreasing
\begin{equation}\label{Ibar-new}
\bar{I}(r)=\frac{1}{r^n}\int_{\bar{E}_r}\left(|\nabla \log
\bar{K}|^2-(\log \bar{K})_\tau\right)\, d\mu\, dt,
\end{equation}
which is just a different appearance of Ecker's quantity \cite{E}.
Remarkably, the integrands in the above monotonic quantities are
again the Li-Yau's expression. Notice that
(\ref{I-new})-(\ref{Ibar-new}) also follow from (\ref{J-new}) and
its cousins by the integrations. In the case of Ricci flow, it was
shown by Perelman in \cite{P} that
\begin{equation}\label{ly-rcf}
|\nabla \log \hat{K}|^2-(\log \hat{K})_\tau=|\nabla
\ell|^2+\ell_\tau
+\frac{n}{2\tau}=\frac{n}{2\tau}-\frac{1}{2\tau^{\frac{3}{2}}}\mathcal{K}(y,
\tau)
\end{equation}
where $\mathcal{K}(y,
\bar{\tau})=\int_0^{\bar{\tau}}\tau^{\frac{3}{2}}H(\frac{d
\gamma}{d\tau})\, d\tau$ with $\gamma$ being the minimizing
$\mathcal{L}$-geodesic joining to $(y, \bar{\tau})$ (assuming that
$(y, \bar{\tau})$ lies out of cut locus), $H(X)=-\frac{\partial
R}{\partial \tau}-\frac{R}{\tau}-2\langle X, \nabla R\rangle
+2Ric(X, X)$. In the case of the mean curvature flow (see \cite{E})
we have that
\begin{equation}\label{ly-mcf}
|\nabla \log \bar{K}|^2-(\log
\bar{K})_\tau=\frac{n}{2\tau}+\langle\vec{H}, \nabla^{\perp} \log
\bar K \rangle-|\nabla^{\perp} \log \bar K|^2.
\end{equation}

When $g(t)$ become singular or degenerate at $(x_0, 0)$, special
cares  are needed in justifying the above identities. Interesting
cases are the gradient shrinking solitons for the Ricci
flow/homothetically shrinking solutions for mean curvature flow
respectively, which have singularity at $\tau=0$. For the case of
the mean curvature flow, one can work in terms of the image of
$y_t(\cdot)$. A homothetically  shrinking solution satisfies that $
\vec{H}(y,\tau)=-\frac{(y-x_0)^{\perp}}{2\tau}$. As in \cite{E}
denote the space time track of $(M_t)$ by
$$
\mathcal{M}=\cup_{t\in (\alpha, \beta)}M_t\times \{t\}\subset
\mathbb{R}^{n+k}\times \mathbb{R}.
$$
It can be checked (cf. \cite{E}) that on a shrinking soliton
$$
\left(\frac{\partial}{\partial \tau}-\Delta +|\vec{H}|^2\right)
\bar{K}(y, \tau; x_0, 0)=0.
$$
For homothetically shrinking solutions, it was shown in \cite{E}
that
$$
\bar{I}(r)=\int_{M_1} \bar{K}(y,1; x_0, 0)\, d\mu_1.
$$
 Apply the space-time divergence theorem to $\widetilde{X}=\nabla
\bar{K}+\bar{K}\frac{\partial}{\partial t}$. We can check that $
\lim_{s\to 0}\int_{P_2^s\cap \mathcal{M}}\, d\mu_s =0 $,  if the
shrinking soliton is properly embedded near$(x_0, 0)$, by
Proposition 3.25 of \cite{E2}. By the virtue of Theorem \ref{mv-sph}
we have that
$$
\bar{J}(r)=\lim_{s\to 0}\int_{P_2^s\cap \mathcal{M}} \bar{K}(y, s;
x_0, 0)\, d\mu_s.
$$
Using the scaling invariance of $\mathcal{M}$ one can show that
$$
\lim_{s\to 0}\int_{P_2^s\cap\mathcal{M}} \bar{K}(y, s; x_0, 0)\,
d\mu_s=\int_{M_1} \bar{K}(y,1; x_0, 0)\, d\mu_1.
$$
Here, with a little abuse of the notation, we denote
$M_1=\mathcal{M}\cap\{\tau=1\}$. Summarizing we have that both
$\bar{I}(r)$ and $\bar{J}(r)$ are  equal to the the so-called
Gaussian density $\Theta(\mathcal{M}, x_0, 0)=\int_{M_t}\bar{K}\,
d\mu_t.$

The Ricci flow case is similar. Assume that $(M ,g(t))$ is a
gradient shrinking soliton (see for example \cite{CLN} for a precise
definition, and \cite{FIK} for the new examples of noncompact
shrinkers). In particular, there exists a smooth function $f(x,
\tau)$ ($\tau=-t$) so that
$$
R_{ij}+f_{ij}-\frac{1}{2\tau} g_{ij}=0.
$$
Moreover, the metric $g(\tau)=\tau \varphi^* (\tau)g(1)$, where
$\varphi(\tau)$ is a one parameter family of diffeomorphisms,
generated by $-\frac{1}{\tau}\nabla_{g(1)}f(x, 1)$. We further
assume that there exists a `attracting' sub-manifold $S$ such that
 $(\nabla f)(x_0, \tau)=0$ (which is equivalent to that $(\nabla
f)(x_0, 1)=0$), for every $x_0\in S$, $f$ is constant on $S$ and the
integral curve $\sigma(\tau)$ of the vector field $\nabla f$ flows
into $S$ as $\tau\to 0$. (All known examples of gradient shrinking
Ricci soliton seem to satisfy the above assumptions). In this case
for any $x_0\in S$ we can define the {\it reduce distance} as
before. However, we require that the competing curves $\gamma(\tau)$
satisfies that $\lim_{\tau\to 0}\gamma'(\tau)\sqrt{\tau}$ exists
(this is the case if $(x_0, 0)$ is a regular time and
$\gamma'(\tau)$ is a $\mathcal{L}$-geodesic. We may define $\ell(y,
\tau)=\inf_{x_0\in S}\ell^{(x_0, 0)}(y, \tau)$. Under this
assumptions, after suitable normalization on $f$ (by adding a
constant), we can show that $f(y, \tau)=\ell(y, \tau)$ (see for
example, \cite{CHI} and
 \cite{CLN}). It can also be easily checked that
 $\hat K$ is a solution to the {\it conjugate heat equation}. Then the
  very similar argument as above shows
 that
$$
\hat{J}(r)=\lim_{s\to 0} \int_{P_2^s} \hat{K}(y, s; x_0, 0)\, d\mu_s
$$
which then implies that $\hat{J}(r)=\int_{M_{1}} \hat{K}\, d\mu_1$,
a constant independent of $r$. Using that relation that
$\hat{I}(r)=\frac{\int_0^r \eta^{n-1}\hat{J}(\eta)\, d\eta}{\int_0^r
\eta^{n-1}\, d\eta}$, we have that $\hat{I}(r)$ also equals to the
`reduced volume' of Perelman $\theta^{(x_0, 0)}(g, \tau)=\int_{M}
\hat{K}\, d\mu_\tau$ (cf. \cite{EKNT}).

\bibliographystyle{amsalpha}

\end{document}